\documentclass[final]{siamltex}
\usepackage{animate}
\usepackage{amsmath}
\usepackage{amssymb}
\usepackage{graphics}
\usepackage{graphicx}
\usepackage{textcomp}
\usepackage{mathrsfs}
\usepackage{epstopdf}
\usepackage{array}
\usepackage{cite}
\usepackage[maxfloats=99]{morefloats}
\usepackage{color}
\usepackage{url}
\usepackage{color}
\usepackage{cite}
\usepackage{subcaption}

\usepackage{stmaryrd}
\usepackage{bm}

\def\e{{\varepsilon}}

\def\T{{\mathcal T}}
\def\N{{\mathcal N}}
\def\L{{\mathcal L}}

\renewcommand{\theequation}{\arabic{section}.\arabic{equation}}

\graphicspath{{fig/}}

\newtheorem{remark}{Remark}[section]

\title{Hessian recovery based  finite element methods  for the Cahn-Hilliard Equation}

\author{ Minqiang Xu\thanks{School of Data and Computer Science, Sun Yat-sen University, Guangzhou 510275, China (xumq9@mail2.sysu.edu.cn).}
\and Hailong Guo\thanks{School of Mathematics and Statistics,  The University of Melbourne,  Parkville, VIC 3010, Australia   (hailong.guo@unimelb.edu.au).}
\and Qingsong Zou\thanks{
School of Data and Computer Science and Guangdong Province Key Laboratory of Computational Science, Sun Yat-sen University, Guangzhou 510275, China
(mcszqs@mail.sysu.edu.cn)}}

\begin{document}

\maketitle

%
%
\medskip

\begin{abstract}
In this paper, we propose a novel recovery based finite element method for the Cahn-Hilliard equation. One distinguishing feature of
the method is that we discretize the fourth-order differential operator in a standard $C^0$  linear finite elements space.
Precisely, we first transform the fourth-order Cahn-Hilliard equation to its variational formulation in which only first-order and second-order derivatives are involved and then
we compute the first  and second-order derivatives of a linear finite element function by a  least-square-fitting recovery procedure.
 When the underlying mesh is uniform of regular pattern, our recovery scheme for the Laplacian operator coincides with  the well-known five-point stencil.  
Another  feature of the method is some  special  treatments on  Neumann type boundary conditions for reducing computational cost. 
 The optimal-order convergence properties and energy stability are numerically proved through a series of benchmark tests.
 The proposed method can be regarded as a combination of the finite difference scheme and the finite element scheme.

\vskip .3cm

\end{abstract}

\section{Introduction}
The phase field model is a powerful tool to characterize interfacial problems in which the dynamics of the physical systems are described by a gradient flow.
The Cahn-Hilliard equation \cite{CH1958} is a famous phase field model introduced by  Cahn and  Hilliard to model the phase separation in binary alloys.  Later on, it is widely used to model multiphase flow\cite{BCB2003, CS2016},  tumor growth\cite{WLFC2008, OHP2010},  and image impainting \cite{BEG2007} etc.

As a  nonlinear parabolic type equation, the analytic solution of the Cahn-Hilliard equation is usually hard to be obtained.  Numerical simulation looks like to be the only feasible way to study the physical problems governed by the Cahn-Hilliard equation. During the past several decades,
a huge number of numerical methods have been developed in the literature, including finite different methods\cite{CR2007}, spectral methods\cite{SY2010}, and finite element methods\cite{WKG2006,ZW2010, EL1992}. In this paper, we concentrate on finite element methods for the Cahn-Hilliard equation.

One of the main difficulties in the numerical solution of the Cahn-Hilliard equation is the discretization of the fourth-order differential operator
in a certain finite element space.
In the  literature,  finite element methods for fourth-order elliptic equations can be roughly categorized into
the following four classes: conforming finite element methods \cite{Ci2002, BS2008}, nonconforming finite element methods \cite{Mo1968},
mixed finite element methods\cite{CR1974}, discontinuous Galerkin method \cite{EGHLMT2002}.  In corresponding,   the Cahn-Hilliard equation has been numerically solved by conforming finite element methods in
\cite{EZ1986, DN1991},  nonconforming finite element methods in \cite{EF1989, ZW2010},  mixed finite element methods in \cite{EFM1989, FP2004}, and  discontinuous Galerkin methods
\cite{WKG2006, XXS2007}.
All the above methods in the primary form
 discretize the Cahn-Hilliard equation at least in a quadratic finite element space,
which means that there are at least six degrees of freedom on each triangular element.
To reduce the complexity,  a new class of finite element methods, called recovery based finite element
methods \cite{CGZZ2017, GZZ2018,  L2014}, is proposed to simulate fourth-order partial differential equations.
The key idea of those methods is to facilitate the simplest element, the continuous linear element, to discretize the fourth-order differential operator.
It is known  that the second-order derivative of $C^0$ piecewise linear function is
not well-defined and thus usually we can not solve a fourth-order differential equation in a $C^0$ linear finite element space.
Such a barrier is alleviated by using the classical gradient recovery operator $G_h$ \cite{ZN2005, ZZ1992a} to smooth the discontinuous piecewise constant function into a continuous piecewise linear function i\cite{CGZZ2017, GZZ2018}.

In this paper, we will also discretize  the
Cahn-Hilliard equation only in the simplest linear element space.
Comparing to the recovery-based FEMs in \cite{CGZZ2017, GZZ2018,  L2014},
the difference here is that we directly recover  the Hessian matrix of a linear finite element
function instead of recovering its gradient.
Note that the Hessian recovery has been studied for the purpose of post-processing \cite{GZZ2017b, AV2002, PABG2011, VMDDG2007}.
In particular, in \cite{GZZ2017b},  Guo et al.  proposed a new Hessian recovery method and established its complete
superconvergence theory on mildly unstructured meshes and ultraconvergence theory on structured meshes.
The Hessian recovery technique in \cite{GZZ2017b} is then applied to solving a sixth-order PDE in  \cite{GZZ2018b}.
In this paper, we use a Hessian recovery technique by firstly recovering a local quadratic polynomial and
then taking second order derivatives of the recovered polynomial as the {\it second-order derivatives} of the linear finite element function.
We sprucely discover that there is an intrinsic connection between the Hessian recovery method and the finite difference method.
In specific, we find that the Hessian recovery method reproduces the standard five-point finite difference scheme on regular pattern uniform meshes.
This means, on the regular pattern uniform meshes, the new recovery based finite element method is a kind of
infusion of the finite difference method and the finite element method in the sense that we first use
the standard five-point finite difference scheme to discretize the Laplacian operator and then put it back
into the standard linear finite element framework.

Different from second-order elliptic equations,  the Neumann boundary condition for fourth-order partial differential equations
is an essential boundary condition in the sense that the boundary condition should be enforced in their associate solution spaces.
But it looks like impossible to be enforced into a $C^0$ linear finite element space.
  In our previous paper \cite{CGZZ2017, GZZ2018, GZZ2018b},  it is imposed by the penalty method \cite{Co1942,ZTZ2013} and the Lagrange multiplier method \cite{Ba1972}.
  Like it for the second order partial differential equations,  the resulting linear system of the penalty method is usually ill-conditioned.
  The Lagrange multiplier method shows the potential to overcome such drawback but it introduces additional degrees of freedom.
  In this paper,  we adopt two different methods to deal with the Neumann boundary condition.
  On general unstructured meshes,   we propose to impose the boundary condition weakly based on a technique called Nitsche's method \cite{Ni1971}, which is originally introduced by Nitsche to incorporate the Dirichlet boundary condition for second order elliptic equations.    A Nitsche's variational formulation for the Cahn-Hilliard equations is presented.
  It paves the way for implicitly imposing the Neumann boundary conditions in Hessian recovery based finite element methods.
As mentioned in the previous paragraph, the Hessian recovery method reduces  to the standard five-point standard finite difference
scheme on uniform meshes.  Such key observation enables us to incorporate the Neumann boundary condition into
the Hessian recovery operator by the celebrated ghost point method in finite difference methods \cite{Le2007}.

The rest of paper is organized as follows: In Section 2, we present a simple introduction to the Cahn-Hilliard
equation and review some of its property.  In Section 3, we first revisit the Hessian recovery method and
uncover its relationship with the classical finite difference method; then, we propose
the new recovery based finite element method to discretize the spatial variable;
the fully discrete formulation is discussed through an energy stable time stepping method.
The proposed method is numerically verified and validated using a series of benchmark examples in Section 4.
We end with some conclusive remarks in Section 5.

\section{The Cahn-Hilliard equation and its variational formulations}
Let $\Omega$ be a bounded polygonal domain with  Lipschitz boundary
$\partial \Omega$ in $\mathbb{R}^2$. For a subdomain $\mathcal{A}$
of $\Omega$, let $\mathbb{P}_m(\mathcal{A})$ be the space of polynomials of
degree less than or equal to $m$ over $\mathcal{A}$ and $n_m$ be the
dimension of $\mathbb{P}_m(\mathcal{A})$ with $n_m=\frac{1}{2}(m+1)(m+2)$.
We denote by $H^k(\mathcal{A})$  the Sobolev space with norm
$\|\cdot\|_{k, \mathcal{A}} $ and seminorm $|\cdot|_{k,  \mathcal{A}}$.
%

The well-known Cahn-Hilliard equation  on a space domain $\Omega$ and a certain time period
$[0,T]$ can be described as below:
\begin{eqnarray}\label{equ:model}
\left \{
\begin{array}{lll}
\frac{\partial u}{\partial t}=-\varepsilon^2 \Delta^2 u+\Delta F'(u),~~~~~~~~&\text{in}~\Omega \times[0,T],\\
\partial_{\mathbf{n}}u=\partial_{\mathbf{n}}(-\varepsilon^2 \Delta u +F'(u))=0,~~~~&\text{on}~\partial \Omega\times[0,T],\\
u(\cdot,0)=u_0(\cdot),~~~~~~~~~~~~~~~~~~&\text{in}~\Omega,\\
\end{array}
\right.
\end{eqnarray}
where $\mathbf{n}$ is the unit outer normal vector  of $\partial \Omega$. The unknown function $u$ often indicates the concentration of one of the two metal components
constituting the alloys, $\e$ is  the size of the interface of two alloys, $F$ is a double well nonconvex function
defined as
$
F(u)=\frac{1}{4}(u^2-1)^2.
$

The Cahn-Hilliard equation \eqref{equ:model} can be viewed as an $H^{-1}$-gradient flow of the  Ginzburg-Landau
free energy functional
\begin{equation}
 E(u):=\int_{\Omega} \left(\dfrac{\varepsilon^2}{2}|\nabla u|^2+F(u)\right) dx,
\end{equation}
of which the first part is called the {\it interfacial energy} and the second part is called the {\it bulk energy}.
Thanks to the homogeneous Neumann boundary conditions, it is easy to verify
that the  mass conservation property
\begin{equation*}
  \frac{d}{dt}\int_\Omega u dx =0,
\end{equation*}
and the energy decay property
\begin{equation*}
~\frac{dE(u)}{dt}=-\|-\varepsilon^2 \Delta u +F'(u))\|^2\leq 0,~\forall~t>0,
\end{equation*}
always hold for the solution of the Cahn-Hilliard equation.

Let $f = F'$. To implicitly impose the Neumann boundary condition $\partial_{\mathbf{n}}u=0$, we introduce  the bilinear form
\begin{equation}\label{equ:bilinear}
a_1(w, v) = \int_{\Omega} \Delta w \Delta vdz
- \int_{\partial \Omega} \Delta w \partial_{\mathbf{n}} v  ds
- \int_{\partial \Omega}   \partial_{\mathbf{n}} w \Delta v ds
+ \gamma\int_{\partial \Omega}\partial_{\mathbf{n}} w  \partial_{\mathbf{n}} v  ds, \forall v,w\in H^2(\Omega),
\end{equation}
with  $\gamma$ is a positive stability parameter to be specified in the
sequel.
It is easy to verify that if $u\in H^2(\Omega)$ is the solution of \eqref{equ:model}, then $u$ satisfies
\begin{equation}\label{equ:nit}
  \left(\frac{\partial u}{\partial t},v\right)+\varepsilon^2a_1(u, v)+(\nabla f(u), \nabla v)=0, \quad  \forall v \in H^2(\Omega).
\end{equation}
Conversely, if $u\in C^4(\Omega)$ and $u_t \in C(\Omega)$ which satisfies \eqref{equ:nit}, then for all $v\in H^2(\Omega)$, we have
\begin{eqnarray}
&&\left(\frac{\partial u}{\partial t}+\varepsilon^2 \Delta^2 u-\Delta f(u),v\right)\nonumber\\
 &&-\int_{\partial \Omega}   \partial_{\mathbf{n}} u \Delta v ds
+ \gamma\int_{\partial \Omega}\ \partial_{\mathbf{n}} u  \partial_{\mathbf{n}} v+\int_{\partial\Omega}\partial_{\mathbf{n}}(-\varepsilon^2 \Delta u +F'(u))v  ds=0.\label{eq:181020}
\end{eqnarray}
Choosing $v\in C_0^\infty(\Omega)$ in \eqref{eq:181020}, we obtain
\[
\left(\frac{\partial u}{\partial t}+\varepsilon^2 \Delta^2 u-\Delta f(u),v\right)=0.
\]
which implies $\frac{\partial u}{\partial t}+\varepsilon^2 \Delta^2 u-\Delta f(u)=0$.
Consequently, \eqref{eq:181020} becomes
\[
-\int_{\partial \Omega}  \partial_{\mathbf{n}} u \Delta v ds
+ \gamma\int_{\partial \Omega}\ \partial_{\mathbf{n}} u  \partial_{\mathbf{n}} v+\int_{\partial\Omega}\partial_{\mathbf{n}}(-\e^2\Delta u+f(u))v  ds=0,v\in H^2(\Omega).
\]
Since $v\in H^2(\Omega)$ is arbitrary, we derive from the above equation that
\[
\partial_{\mathbf{n}} u=\partial_{\mathbf{n}} (-\e^2\Delta u+f(u))=0
\]
It means  $u$ is the classical solution of  \eqref{equ:model}. From the above reasonings, we obtain that
the solution of \eqref{equ:nit} is a weak solution of \eqref{equ:model} and we call \eqref{equ:nit}
a variational formulation of \eqref{equ:model}.

\begin{remark}An alternative variational equation of \eqref{equ:model} is
\begin{equation}\label{equ:nit2}
  (\frac{\partial u}{\partial t},v)+\varepsilon^2a_2(u, v)+(\nabla f(u), \nabla v)=0, \quad  \forall v \in H^2(\Omega)
\end{equation}
where \begin{equation}\label{equ:bilhes}
   a_2(w, v) = \int_{\Omega} D^2w: D^2 vdx
- \int_{\partial \Omega} \partial^2_{\mathbf{n}} w\partial_\mathbf{n} v ds
- \int_{\partial \Omega}  \partial_\mathbf{n} w\partial^2_\mathbf{n} v  ds
+ \gamma \int_{\partial \Omega}\ \partial_\mathbf{n} w \partial_\mathbf{n}v ds,
\end{equation}
and $ A:B$ is the Frobenius norm of $2\times 2$ matrices.
We observe that here, the bilinear form $a_2(\cdot,\cdot)$ differs from $a_1(\cdot,\cdot)$ by replacing the Laplace operator $\Delta$ with the
Hessian matrix operator $D^2$.

In both variational formulations, the Neumann boundary conditions $\partial_{\mathbf{n}}u=0$ is weakly built into the bilinear formations.
The idea is similar to impose the Dirichlet boundary condition for the second-oder elliptic equations  by Nitsche \cite{Ni1971}.
We call those two methods  the Nitsche's method.
\end{remark}

\begin{remark}
Originally,  the Neumann  boundary condition for fourth-order partial differential equations is enforced into their solution space.  For such purpose, let
\begin{equation}\label{equ:subspace}
 V = \{ v\in H^2(\Omega):  \partial_{\mathbf{n}}v = 0 \text{ on } \partial \Omega\}.
\end{equation}
The bilinear forms $a_1$ and $a_2$ in $V$ reduce to
\begin{equation}\label{equ:bilinear3}
a_3(w, v) = \int_{\Omega} \Delta w \Delta vdx, w,v\in V
\end{equation}
and
\begin{equation}\label{equ:bilinear4}
   a_4(w, v) = \int_{\Omega} D^2w : D^2 vdx,w,v\in V
\end{equation}
respectively. Correspondingly, the variational formulations  become to find  $u\in (L^2([0,T];V)$ such that
\begin{equation}\label{equ:var}
  \left(\frac{\partial u}{\partial t},v\right)+\varepsilon^2a_3(u, v)+(\nabla f(u), \nabla v)=0, \quad  \forall v \in V,
\end{equation}
or, to find  $u\in (L^2([0,T];V)$ such that
\begin{equation}\label{equ:var2}
  \left(\frac{\partial u}{\partial t},v\right)+\varepsilon^2a_4(u, v)+(\nabla f(u), \nabla v)=0, \quad  \forall v \in V.
\end{equation}
\end{remark}

\section{A novel recovery based finite element method}
In this section, we design novel recovery-technique-based finite element methods for Cahn-Hilliard equations.
Since our main attention  is on  novel space discretization techniques, for the time
discretization, we choose a simple one-step energy stable linear scheme proposed in \cite{ZW2010,SY2010,HLT2007}. Precisely, we make use of a semi-implicit scheme with
an extra stabilized penalty term added to ensure energy stability. Let the time step size be $\Delta t = \frac{T}{N}$,
$u^0(x)= u(x,0)$, and $u^n(x)\approx u(x,n \Delta t), n=1,2,\ldots,N$, then the stabilized first-order semi-implicit
method reads as: find $u^n\in S, n=1,2,\ldots N$ such that for all $v\in S$,
\begin{equation}\label{semiformulation}
\begin{split}
&\left(\frac{u^{n+1}-u^n}{\Delta t},v\right)+\varepsilon^2a_i(u^{n+1}, v) \\
&+(\nabla f(u^n), \nabla v) +\kappa\left(\nabla (u^{n+1}- u^{n}),\nabla v\right)=0, \quad i=1,2,3,4,
\end{split}
\end{equation}
where $S= H^2(\Omega)$ for $i=1,2$ and  $S= V$ for $i=3,4$.
Note that the choice of $\kappa$ has a great influence on the stability of \eqref{semiformulation}, see
 \cite{ZW2010,SY2010, LQT2016, LQ2017} for the details. In this paper, we always take $\kappa=2$.
 In fact,  \cite{LQT2016, LQ2017} give you rigorous analysis on the choose of $\kappa$.

Next we explain how to discretize \eqref{semiformulation} in a linear finite element space. Let $\mathcal{T}_h$ be a shape regular {\it triangulation} of $\Omega$
 with mesh size  $h$.  The set of all vertices  and  of all edges of $\mathcal{T}_h$ are denoted by $\mathcal{N}_h$ and $\mathcal{E}_h$, respectively.
We define the  standard continuous linear finite element space $S_h$ on $\mathcal{T}_h$  by
\begin{equation}\label{equ:p1fes}
S_h := \left\{v_h \in C^{0}(\Omega):   v_h|_T \in {\mathbb P}_1,\forall T\in \T_h\right\}.
\end{equation}
and denote its nodal basis by $\{\phi_z\}_{z\in \mathcal{N}_h}$.

To construct our fully discrete schemes on the linear finite element space $S_h$, we first introduce a Hessian recovery technique based on least-squares fitting in
the first subsection. Then we apply the Hessian recovery operator to develop our novel fully
discrete schemes for the Cahn-Hilliard equation in the second subsection.
\subsection{A Hessian recovery operator in linear finite element spaces}
It is known that the second order derivative of a function $v_h\in S_h$ equals to $0$ in the interior
of each element $T\in\T_h$ and is not well-defined on an edge $E\in {\mathcal E}_h$. In the following, we propose a
{\it least-square type} method to calculate the {\it approximate second order derivatives} of $v_h$. In other words,
we will define a  {\it Hessian recovery operator} $H_h$  from $S_h$ to $S_h^4$ which maps a function $v_h\in S_h$ to $H_hv_h\in S_h^4$
so that $H_hv_h$ can be regarded as an approximation of the Hessian matrix of $v_h$ in some sense.

Since $H_hv_h\in S_h^4$, to define $H_hv_h$, it is sufficient to  define the value  $(H_hv_h)(z)$ for all $z\in \mathcal{N}_h$.
For this purpose,  we first construct a local {\it patch} associated with $z$ which is a polygon surrounding the node $z$.
Given a vertex $z\in {\mathcal N}_h$ and a nonnegative integer $n\in \mathbb{N}$,  let the first $n$ layer element patch be
\begin{equation}
\mathcal{L}(z,n)=
\begin{cases}
\{z\}, & \text{if}\ n=0,\\
\bigcup\{T: T\in\T_h,\ \overline{T}\cap\mathcal{L}(z,0)\neq\phi\}, &
\text{if}\ n=1,\\
\bigcup\{T: T\in\mathcal{T}_h,\ \overline{T}\cap\overline{\mathcal{L}(z,n-1)}\
\text{is an edge in } \mathcal{E}_h\}, & \text{if}\ n\ge 2.
\end{cases}
 \label{local}
\end{equation}
For all $z\in\N_h$, let $n_z$ be the smallest integer  such that $\mathcal{L}(z,n)$ satisfies the rank condition
in the following sense.
\begin{definition}\label{def:rank}
 A surrounding $z$ polygon  is said to satisfy the {\it rank condition} if it admits
 a unique least-squares fitted polynomial $p_{z}$ in \eqref{equ:pprls}.

\end{definition}

We define the local patch associated with $z$ as $\Omega_z=\L(z,n_z)$.
Using the vertices in $\Omega_{z}$ as
sampling points, we fit a quadratic polynomial  $p_{z}$ at the vertex $z$ in the following
least-squares sense
 \begin{equation}\label{equ:pprls}
p_{z}=\arg\min_{p\in \mathbb{P}^2(\Omega_{z})} \sum_{x\in \in \overline{\Omega}_{z}\cap \mathcal{N}_h}| {p(x)-v_{h}(x)}|^2.
\end{equation}
Then, we define the recovered Hessian node value by
\begin{equation}\label{equ:hessian}
(H_hv_h)(z)
=\left(
\begin{matrix}
 H_h^{xx}v_h(z) &  H_h^{xy}v_h(z)\\
  H_h^{yx}v_h(z) &  H_h^{yy}v_h(z)
\end{matrix}
\right)=\left(
\begin{matrix}
 \frac{\partial^2 p_{z}}{\partial x^2}(z) & \frac{\partial^2 p_{z}}{\partial x\partial y}(z) \\
 \frac{\partial^2 p_{z}}{\partial y\partial x}(z)&  \frac{\partial^2 p_{z}}{\partial y^2}(z)  \\
\end{matrix}
\right).
\end{equation}
With this definition, we have $H_hv_h=\sum_{z\in\N_h} H_hv_h(z)\phi_z$ and the symmetric property $H_h^{xy}=H_h^{yx}$
of the Hessian matrix function $H_hv_h$.
Moreover, based on $H_h$,  we  define  a discrete Laplacian operator $\Delta_h: S_h \rightarrow S_h$ as
\begin{equation}\label{equ:lap}
\Delta_hv_h  = H_h^{xx}v_h + H_h^{yy}v_h.
\end{equation}
Note that in the same way, we can recover the {\it gradient} of $v_h$ by letting
\begin{equation}
(G_h v_h)(z) = \nabla p_{z}(z),\forall z\in\N_h
\end{equation}
and $G_hv_h=\sum_{z\in\N_h} G_hv_h(z)\phi_z\in S_h^2.$

\begin{figure}[!h]
   \centering
   \subcaptionbox{\label{fig:regular}}
  {\includegraphics[width=0.4\textwidth]{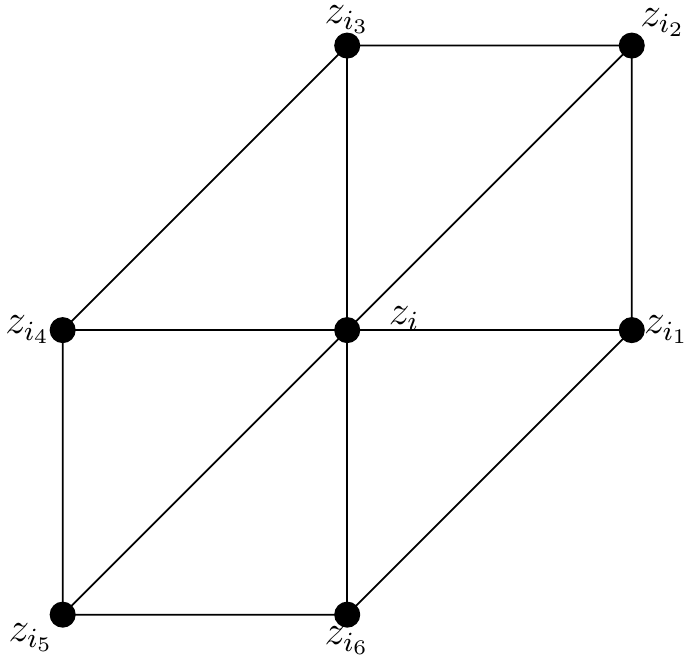}}
  \subcaptionbox{\label{fig:laplace}}
   {\includegraphics[width=0.4\textwidth]{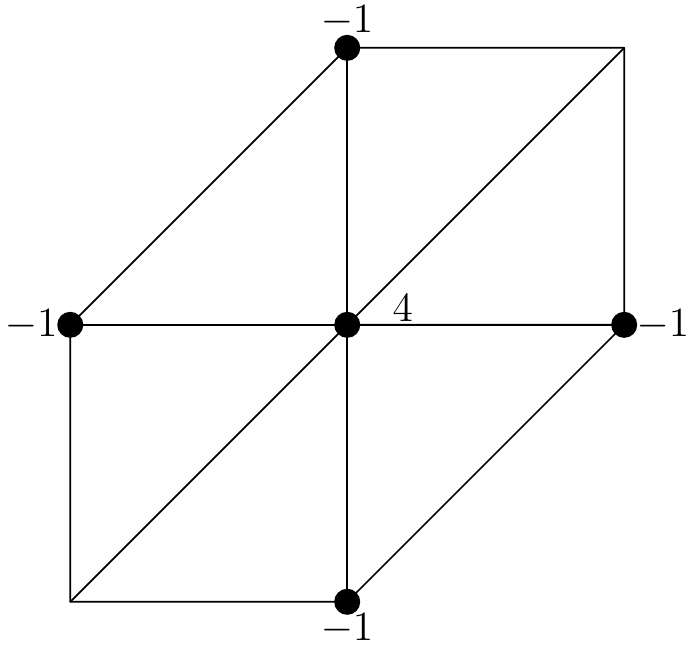}}
   \caption{Illustration of Hessian recovery on uniform mesh: (a) local patch ;  (b) Discrete Laplace operator.}
\end{figure}

Note that the Hessian recovery operator $H_h$ in \eqref{equ:hessian} has been applied to post-process the finite element solution in \cite{GZZ2017b}.
 According to the numerical results, sometimes it might be inefficient or even not convergent as a post-processing operator.
 Since it involves a relatively smaller number of neighbourhood vertices in its stencil, here we choose it as our pre-processing operator.
In fact, $H_h$ can be regarded as a special finite difference operator of the second order on uniform meshes and of the first order on general unstructured meshes.
To elucidate this basic idea, we consider a special case that $\T_h$ is a regular pattern uniform triangular mesh. In this case, for an interior node $z=z_i$, the local patch
$\Omega_{z_i}$ is defined as the polygon $z_1\cdots z_6$, see  Figure \ref{fig:regular} and
thus the sampling points in $\Omega_{z_i}$ include $z_{i_0}, z_{i_1}, \cdots, z_{i_6}$ with $z_{i_0}=z_{i}$.  Using these seven sampling points, we
 fit a quadratic polynomial $p_{z_i}$ in the least-squares sense and take second-order differentiation,
 which produces
 \begin{align*}
  (H_h^{xx}u)(z_i) = &\frac{1}{h^2}(u_1-2u_0+u_4),\\
  (H_h^{xy}u)(z_i) = &\frac{1}{2h^2}(2u_0-u_1+u_2-u_3-u_4+u_5-u_6),\\
    (H_h^{xx}u)(z_i) = &\frac{1}{h^2}(u_3-2u_0+u_6);
\end{align*}
 where the values $u_j$ is defined by $u_j = u(z_{i_j})$.   By the definition of the discrete
 Laplacian  operator \eqref{equ:lap},   we have
 \begin{equation}\label{equ:reglap}
(\Delta_h u)(z_i) = \frac{1}{h^2}(u_1+u_3-4u_0+u_4+u_6).
\end{equation}
The formula \eqref{equ:reglap} implies the discrete Laplacian operator on regular pattern uniform meshes
is the well-known five-point-finite-difference stencil of the Laplace operator,  as illustrated in Figure \ref{fig:laplace}.
By the standard approximation theory, we have
 \begin{equation}\label{equ:regerror}
\| \Delta u - \Delta_h u_I\|_{0, \Omega} \le C h^2\|u\|_{4, \Omega}.
\end{equation}
This exciting discovery implies that the Hessian recovery operator can be regarded as an extension
of the classic second-order difference operator on regular meshes to the
difference operator on non-uniform meshes and thus it can be used to design discrete schemes
for higher-order differential equations.

\begin{remark}\label{rmk:bndpatch}
 For a boundary vertex $z\in \partial \Omega$, there are other  approaches to construct the local patch $\Omega_z$, see \cite{GZZZ2016} for the details.
\end{remark}

\subsection{Fully discrete schemes}
%

%
To present our fully discrete schemes, we first introduce the discrete bilinear form $a_{i,h}(\cdot, \cdot)$ on the linear finite element space $S_h$.
 For all $w_h, v_h \in S_h$, we define
\begin{eqnarray}
 a_{1,h}(w_h, v_h) &= & \int_{\Omega} \Delta_h w_h \Delta_h v_hdz
- \int_{\partial \Omega} \Delta_h w_h (G_hv_h\cdot \mathbf{n})  ds\nonumber\\
&&-\int_{\partial \Omega}(G_hw_h\cdot \mathbf{n})\Delta_h v_h ds
+\gamma\int_{\partial \Omega}(G_hw_h\cdot \mathbf{n})(G_hv_h\cdot \mathbf{n})  ds,\label{equ:recbilinear}
\end{eqnarray}
and
 \begin{eqnarray}
   a_{2,h}(w_h, v_h) &= &\int_{\Omega} H_hw_h H_hv_hvdz
- \int_{\partial \Omega} (\mathbf{n}^TH_hw_h\mathbf{n}) (G_hv_h\cdot \mathbf{n}) ds\nonumber\\
-& \int_{\partial \Omega}&  (G_hw_h\cdot \mathbf{n})  (\mathbf{n}^TH_hv_h\mathbf{n}) ds
+ \gamma\int_{\partial \Omega}(G_hw_h\cdot \mathbf{n})(G_hv_h\cdot \mathbf{n})  ds,\label{equ:recbilhes}
\end{eqnarray}
where $\gamma= \frac{C}{h}$ with  $C$  a sufficiently large positive constant.

The fully discrete Hessian recovery based finite element method for the Cahn-Hillliard equation \eqref{equ:model}
 reads as :
find $\{u_h^n\}_{n\geq 1} \in V_{h}$ such that for all  $v_h \in S_h$,
\begin{equation}\label{equ:recnit}
 \left(\frac{u^{n+1}_h-u^n_h}{\Delta t},v_h\right)+\varepsilon^2a_{i,h}(u_h^{n+1}, v_h)+\kappa (\nabla u_h^{n+1}- \nabla u_h^n,\nabla v_h)
 +(\nabla f(u_h^n), \nabla v_h)=0.
\end{equation}
%
%

Note that both the schemes in \eqref{equ:recnit} work on general unstructured meshes. Moreover,
our numerical experiments indicate that there is no essential difference between these two schemes.
We also observe that the stiffness matrices corresponding to  both schemes \eqref{equ:recnit}
are symmetric and positive definite, so both schemes are stable and uniquely solvable.

%
%
Next,  we present a simple fully discrete scheme derived from the variational formulation \eqref{equ:var} on uniform meshes.
 We define the finite element space $S_{h,0}\subset S_h$ as
\begin{equation}\label{equ:femsp}
S_{h, 0}  = \{ v_h \in S_h:  \nabla_h v_h(z) \cdot \mathbf{n}=0, \forall z \in \mathcal{N}_h\cap\partial\Omega \},
\end{equation}
where $\nabla_h$ is a discrete gradient operator so that  $S_{h,0}$ be a discrete analogous of $V$.  Note that in the continuous linear finite element space $S_h$,
$\nabla v_h \cdot \mathbf{n}$ is not well defined on a vertex of $\T_h$, so at each boundary vertex, we
use  a central finite difference scheme instead of  $\nabla v_h(z)$ to define $\nabla_h v_h(z)$.

The key part  in construction of a simpler scheme on the uniform  meshes is  based on the fact the discrete Laplacian operator $\Delta_h$ reduces to
the five-point-finite-difference stencil at an interior vertex, as illustrated  in \eqref{equ:reglap}.
Now,  we suppose $\Delta_h$ is the discrete Laplacian operator on the finite element space $S_{h,0}$.
Different from the general treatment of the recovery on the boundary as introduced in the previous subsection,
we borrow the idea of ghost point method from the finite difference method \cite{Le2007}.
In specific, at every boundary vertex $z_i$, we introduce one or more ghost points.
Then, we still have the fact that the  the discrete Laplacian operator $\Delta_h$ is just the five-point finite difference
stencil  at the boundary vertex $z_i$ but it involves the value of the finite element function at the ghost points.
To eliminate it, we combine the discrete boundary condition $\nabla_h v_h(z_i) \cdot \mathbf{n}=0$ which also
involves the same ghost points.

\begin{figure}[!h]
   \centering
   \subcaptionbox{\label{fig:ghost}}
  {\includegraphics[width=0.35\textwidth]{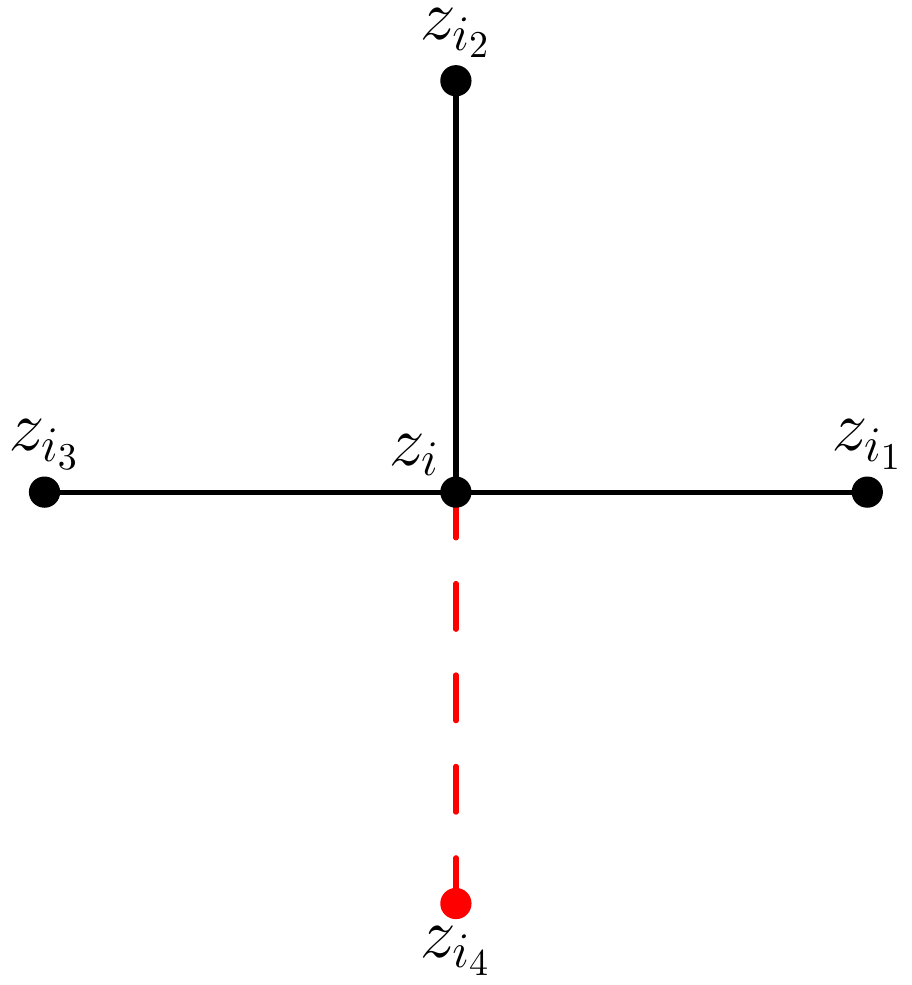}}
  \subcaptionbox{\label{fig:corner}}
   {\includegraphics[width=0.35\textwidth]{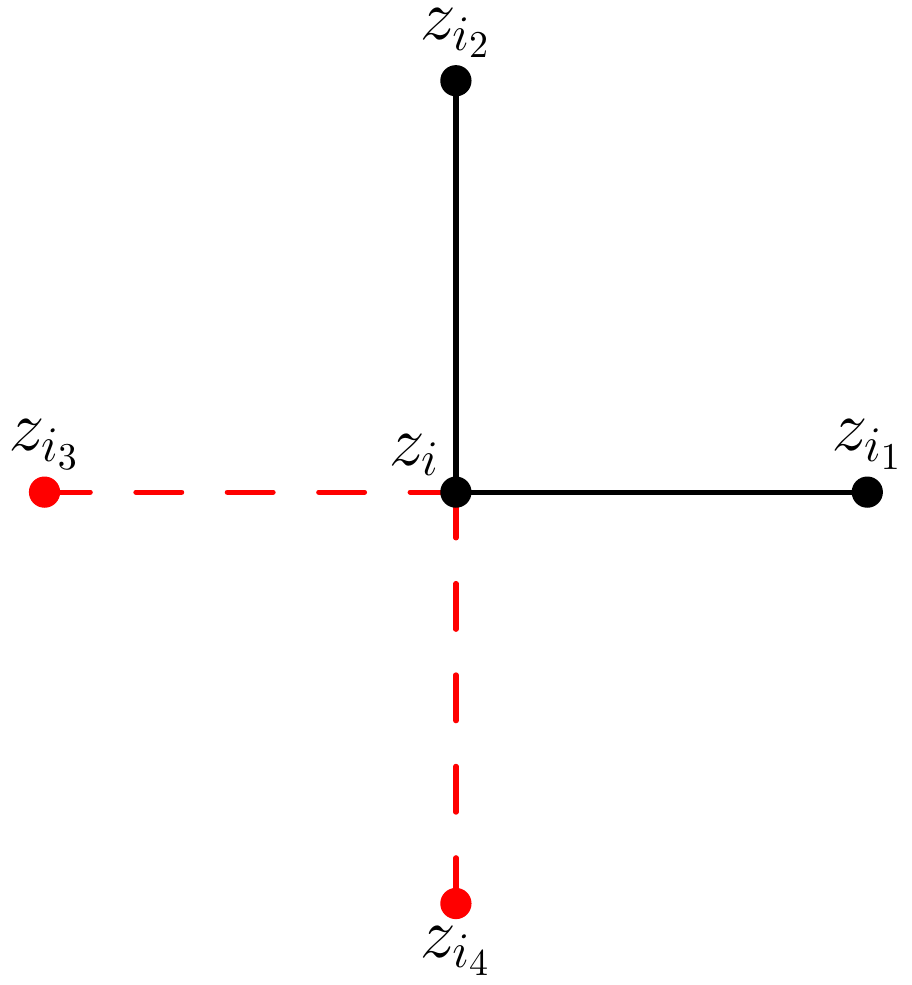}}
   \caption{Illustration of ghost point method: (a). One ghost point; (b) Two Ghost points.}
\end{figure}

To illustrate idea, we consider a typical boundary vertex $z_i$, as illustrated in Figure \ref{fig:ghost}.
In that case, $z_i$ is a boundary vertex with three neighbour mesh vertices $z_{i_1}, z_{i_2}, z_{i_3}$.
To apply the five-point finite difference scheme at $z_i$, we introduce  a ghost point $z_{i_4}$, as
the red dot  point in Figure  \ref{fig:ghost},  and  the discrete Laplacian  $\Delta v_h (z_i) $  is
\begin{equation}\label{equ:bndlap}
\Delta v_h (z_i)   =  \frac{1}{h^2}\left(v_h(z_{i_1})+v_h(z_{i_2})-4v_h(z_i)+v_h(z_{i_3})+v_h(z_{i_4})\right),
\end{equation}
which involves the ghost point finite element function value $v_h(z_{i_4})$.
By the definition of the finite element space $S_{h,0}$,   at the boundary vertex $z_i$, we also
\begin{equation}\label{equ:bndneu}
 \nabla_h v_h(z_i) \cdot \mathbf{n} =   \frac{1}{2h}\left(v_h(z_{i_4})-v_h(z_{i_2})\right) = 0.
\end{equation}
Using \eqref{equ:bndlap} and \eqref{equ:bndneu} to eliminate $v_h(z_{i_4})$, we obtain
\begin{equation}\label{equ:newbndlap}
\Delta v_h (z_i)   =  \frac{1}{h^2}\left(v_h(z_{i_1})+2v_h(z_{i_2})-4v_h(z_i)+v_h(z_{i_3}))\right),
\end{equation}
which only depends on the value at the vertices in $\mathcal{N}_h$.   In other word,
we have embedding the discrete Neumann boundary condition into
the discrete Laplacian operator $\Delta_h$.
Similarly, we can explicitly construct the discrete Laplacian operator $\Delta_h$ at a corner boundary vertex,
which may need  two ghost points as plotted in Figure \ref{fig:corner}.
In this case, the computation of $\Delta_h$ does not need to use an implicit least-squares fitting process.

\begin{remark}
The key processing is to build  the discrete Neumann boundary condition into  the discrete Laplacian operator $\Delta_h$.
Such process  is only possible for the uniform meshes.  For  general unstructured meshes, we can define a similar discrete
finite element space as $S_{h,0} = \{ v_h \in S_h:  G_hv_h(z_i) \cdot \mathbf{n}=0, \forall z_i \in \mathcal{N}_h\cap\partial\Omega \}$.
But the discrete boundary condition $G_hv_h(z_i) \cdot \mathbf{n}=0$ can not be embedded into discrete Laplacian operator  $\Delta_h$.
We have to use other methods like the penalty method \cite{Co1942,ZTZ2013} and the Lagrange multiplier method \cite{Ba1972}  to impose the discrete Neumann boundary condition
$G_hv_h(z_i) \cdot \mathbf{n}=0, \forall z_i \in \mathcal{N}_h\cap\partial\Omega$.  However,  these two methods perform badly for the Cahn-Hilliard equation.
\end{remark}

Then the discrete bilinear $a_{3,h}(\cdot, \cdot)$ on $S_{h,0}$ as
\begin{equation}\label{equ:unibilinear}
\begin{split}
 a_{3,h}(w_h, v_h) = & \int_{\Omega} \Delta_h w_h \Delta_h v_hdx, \quad \forall v_h,w_h\in S_{h,0}.
\end{split}
\end{equation}
A simple fully discrete  for  \eqref{equ:model} on uniform meshes is to find
$\{u_h^n\}_{n\geq 1} \in S_{h,0}$ such that for all $v_h \in S_{h,0}$,
 \begin{equation}\label{equ:simplescheme}
 \left(\frac{u^{n+1}_h-u^n_h}{\Delta t},v_h\right)+\varepsilon^2 a_{3,h}(u_h^{n+1}, v_h)+\kappa (\nabla u_h^{n+1}- \nabla u_h^n,\nabla v_h)\\
+(\nabla f(u_h^n), \nabla v_h)=0.
\end{equation}

Since the bilinear form $a_{3,h}(\cdot,\cdot)$ does not involve the computation of the gradient recovery operator $G_h$ either,
the scheme \eqref{equ:simplescheme}  is very computationally efficient and accurate.
Moreover, the scheme \eqref{equ:simplescheme} can be regarded as a mixture of the finite difference method and the finite element method
since we first use the finite difference operator $\Delta_h$ to recover the second order derivatives of a linear finite element function
and then bring them back to the framework of the finite element method. Namely, the scheme \eqref{equ:simplescheme}
sheds some light on using the finite difference operators to construct simple finite element methods for higher-order partial differential equations.

It may worth mentioning that the gradient recovered method proposed in \cite{GZZ2018}
for fourth-order problems use the minimum number of degrees of freedom among the finite element spaces, see  for details.
However, our numerical experiments show that a direct application of the gradient recovered method to the Cahn-Hilliard equation leads to an unstable
scheme. Moreover, compared with the gradient recovered method,
the present Hessian recovered method uses the same number of total degree
but its stiffness matrix is more sparse than the one  derived from the gradient-recovered method.

\section{Numerical Experiments}
In this section,  we present several numerical examples to demonstrate the properties of our proposed methods.

\begin{figure}[!h]
   \centering
   \subcaptionbox{\label{fig:uniform}}
  {\includegraphics[width=0.3\textwidth]{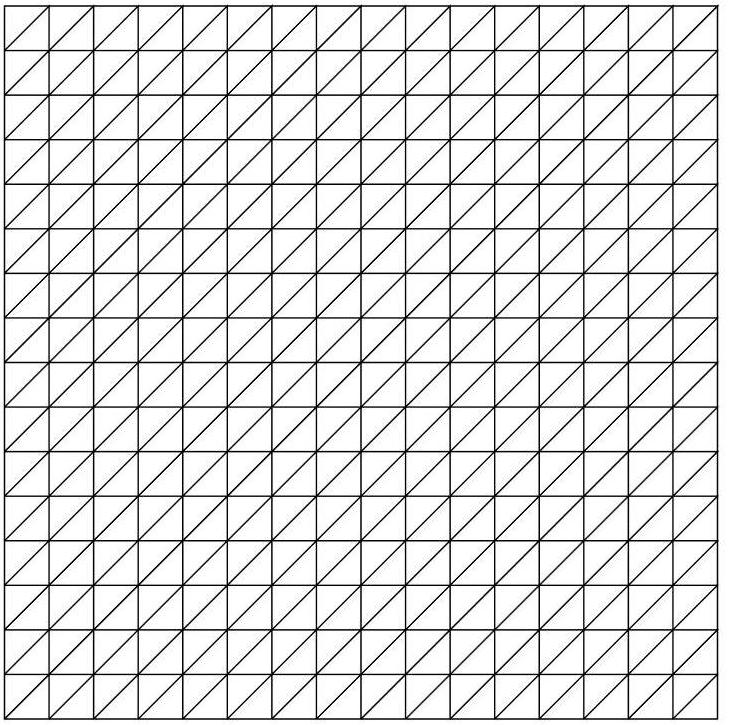}}
  \subcaptionbox{\label{fig:unstruct}}
   {\includegraphics[width=0.3\textwidth]{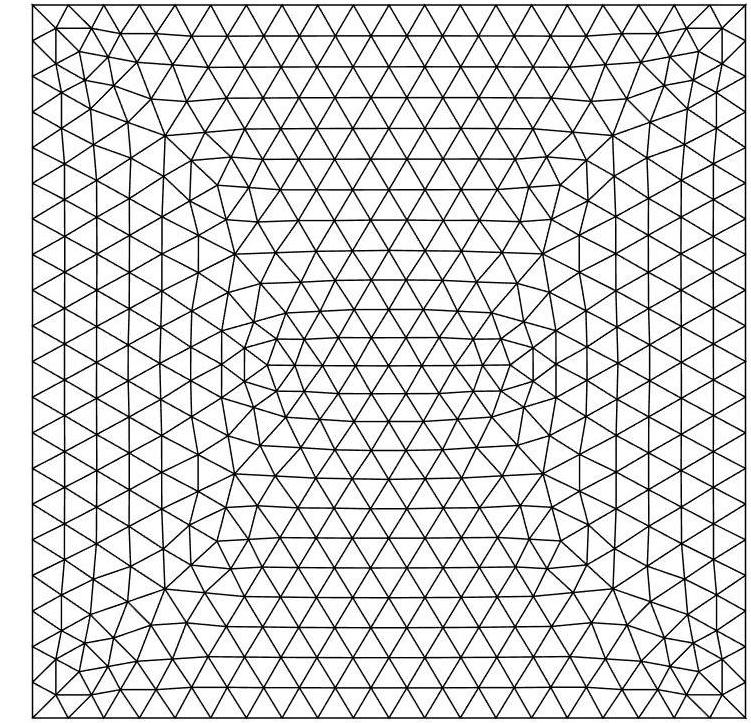}}
   \caption{ (a) A uniform mesh; (b) An unstructured mesh.}
   \label{fig:msh}
\end{figure}

Except for the last numerical example, the domain $\Omega$ of the problems in this section is chosen as the unit square $[0,1]^2$.
In our experiments, we will adopt two different types of meshes:  the uniform and unstructured meshes.
Our uniform meshes are generated by first dividing $\Omega$ into $m^2$ congruent subsquares and then splitting each subsquare into two right-angled triangles, see Figure \ref{fig:uniform}.
Our unstructured meshes are generated by the first partition of  the domain with the Delaunay mesh generator {\it EasyMesh} \cite{easymesh} to obtain the first level mesh and then uniformly refines each triangle in the first level mesh several times, see Figure \ref{fig:unstruct}.
On a uniform mesh, we will use the fully discrete scheme  \eqref{equ:simplescheme}, while on an unstructured mesh, we will use
the scheme \eqref{equ:recnit} with $i=1,\kappa=2,C=1$ to solve numerically the Cahn-Hilliard equations.

Throughout this section,we define  discrete interface energy and discrete buck energy at time $t_n$ respectively as :
\begin{equation*}
E_1^n=\frac{\epsilon^2}{2}\sum_{\tau\in \mathcal{T}_h}\int_\tau |G_h u_h^n|^2\text{d}x\text{d}y,\quad E_2^n=\frac{1}{4}\sum_{\tau\in \mathcal{T}_h}\int_\tau ((u_h^n)^2-1)^2\text{d}x\text{d}y.
\end{equation*}
Moreover, we denote different kinds of numerical errors by
\begin{eqnarray*}
~e_0 &=& \|u-u_h\|_{0,\Omega},
  ~e_1=\|\nabla u-\nabla u_h\|_{0,\Omega}, \\
 ~e _{1,r}&=&\|\nabla u-G_h u_h\|_{0,\Omega},
 e_2=\|D^2 u-H_hu_h\|_{0,\Omega},
\end{eqnarray*}
and use
$r=\frac{\log(e_h/e_{\frac{h}{2}})}{\log(2)}$
to indicate the convergence rates.

\subsection{Accuracy}

{\bf Example 1:} We consider the non-homogeneous Cahn-Hilliard equation
\begin{eqnarray}\label{equ:ex1}
\left \{
\begin{array}{lll}
\frac{\partial u}{\partial t}=-\varepsilon^2 \Delta^2 u+\Delta(u^3-u)+g,~~~~~~~~&\text{in}~\Omega \times[0,T],\\
\partial_{\mathbf{n}}u=\partial_{\mathbf{n}}\Delta u =0,~~~~&\text{on}~\partial \Omega\times[0,T],\\
\end{array}
\right.
\end{eqnarray}
with  the parameter $\epsilon=0.1$. The initial solution $u_0$ and $g$ are chosen such that the exact solution is
\begin{equation*}
  u(x,y,t)=e^{-2t} \cos(\pi x) \cos(\pi y).
\end{equation*}

To compute the convergence rates with respect to the space meshsize $h$, we  fix the time step size $\Delta t=10^{-6}$
and study the convergence order of the numerical solution at  $T=0.1$ computed by  the  scheme \eqref{equ:simplescheme} on uniform meshes and
 by  the scheme (\ref{equ:recnit}) $(i=1,C=1,\kappa=2)$
on  unstructured meshes. The numerical results by \eqref{equ:simplescheme} and
\eqref{equ:recnit} are depicted in  Table 4.1 and Table 4.2, respectively.
From these two tables, we observe that for both schemes, the  $L^2$-norm errors converge with order 2 while the $H^1$-seminorm errors
converge with order $1$ which are both optimal for a linear finite element method.
We also observe that the recovered $H^1$-seminorm error is superconvergent of  $\mathcal{O}(h^2)$ while
 the recovered $H^2$-seminorm errors  converges optimally with order 1.

To test the convergence rate of the scheme (\ref{equ:simplescheme}) with respect to the time discretization, we fix
the spacial mesh size $h=1/128$. The corresponding numerical results at $T=0.01$ with
different time step $\Delta t$ are shown in Table~4.3. The numerical
results evidently indicate that the scheme (\ref{equ:simplescheme})
is of first order in time, which is consistent with the first-order
semi-implicit  scheme.
\begin{table}[!htp]
\vskip -0.5cm
\small{\caption{\emph{Spatial errors and convergence rates by scheme (\ref{equ:simplescheme}) for Example 1}}}
\begin{center}
{\small
\begin{tabular}{l|ll|ll|ll|ll}
\hline
$~~~h $  &$\quad\quad e_0$     &~$r$ &$\quad \quad e_1$   &~$r$  &$\quad\quad e_{1,r}$    &~$r$  &$\quad\quad e_2$ &~$r$ \\
\hline
 $~1/16$ &$ 1.87\times10^{-2}$ &     &$2.46\times10^{-1}$ &      &$1.49\times10^{-1}$ &   &$2.78\times10^{-0}$ &\\
 $~1/32$ &$ 4.09\times10^{-3}$ &2.2  &$9.65\times10^{-1}$ &$1.4$ &$2.99\times10^{-2}$ &2.3&$6.52\times10^{-1}$ &1.9\\
 $~1/64$ &$ 9.92\times10^{-4}$ &2.0  &$4.55\times10^{-2}$ &$1.1$ &$7.47\times10^{-3}$ &2.0&$2.98\times10^{-1}$ &1.3\\
 $1/128$ &$ 2.47\times10^{-4}$ &2.0  &$2.24\times10^{-2}$ &$1.0$ &$1.87\times10^{-3}$ &2.0&$1.39\times10^{-1}$ &1.1\\
 $1/256$ &$ 6.14\times10^{-5}$ &2.0  &$1.12\times10^{-2}$ &$1.0$ &$4.67\times10^{-4}$ &2.0&$6.86\times10^{-2}$ &1.0\\
\hline
\end{tabular}}
\end{center}
\label{tab:10}
\end{table}


\begin{table}[!htp]
\vskip -0.5cm
\small{\caption{\emph{Spatial errors and convergence rates by scheme  (\ref{equ:recnit}) for Example 1}}}
\begin{center}
{\small
\begin{tabular}{l|ll|ll|ll|ll}
\hline
~~\text{dof} &$~\quad e_0$   &~$r$&$~~~\quad e_1$      &~$r$   &$~~~\quad e_{1,r}$  &~$r$  &$\quad \quad e_2$   &~$r$ \\
\hline
~~513  &$1.20\times10^{-2}$  &    &$1.61\times10^{-1}$ &       &$8.77\times10^{-2}$ &      &$2.31\times10^{-0}$ &\\
~1969  &$2.45\times10^{-3}$  &2.3 &$5.64\times10^{-2}$ &$1.5$  &$1.66\times10^{-2}$ &$2.4$ &$1.10\times10^{-1}$ &1.1\\
~7713  &$6.01\times10^{-4}$  &2.0 &$2.59\times10^{-2}$ &$1.1$  &$4.16\times10^{-3}$ &$2.0$ &$5.37\times10^{-1}$ &1.0\\
30529  &$1.52\times10^{-4}$  &2.0 &$1.27\times10^{-2}$ &$1.0$  &$1.08\times10^{-3}$ &$2.0$ &$2.66\times10^{-1}$ &1.0\\
\hline
\end{tabular}}
\end{center}
\label{tab:5}
\end{table}
\begin{table}[!htp]
\small{\caption{\emph{Temporal errors and convergence rate by scheme (\ref{equ:simplescheme}) for Example 1}}}
\begin{center}
{\small
\begin{tabular}{l|llll|l}
\hline
$\Delta t$ &$\quad\quad 10^{-3}$      &$\quad 10^{-3}/2$  &$\quad 10^{-3}/2^2$    &$\quad 10^{-3}/2^3$   &~$r$ \\
\hline
 $e_0$ &$ 1.13\times10^{-3}$ &$5.66\times10^{-4}$       &$2.83\times10^{-4}$ &  $1.41\times10^{-4}$    &1.0\\
\hline
\end{tabular}}
\end{center}
\label{tab:2}
\end{table}

{\bf Example 2:} We consider the Cahn-Hilliard equation (\ref{equ:model})
with the parameter $\epsilon=0.1$ and the initial value $u_0=\cos(\pi x)\cos(\pi y)$.

We use the simple scheme (\ref{equ:simplescheme}) to compute the numerical results.
As the exact solution of Example 2 is unknown, we use the computable quantity $u_{\frac{h}{2}}-u_h$
to replace the  ``true error" $e=u-u_h$ in our real computations.
As in the previous example, we fix $\Delta t=10^{-5}$ and   $T=0.1$
to test the convergence behaviour of the spatial discretization.
The corresponding numerical errors and convergence rates are shown in Table 4.4.
We can observe similar convergence and superconvergence results as in Example 1.

Also as in the previous example, we  test the convergence rate of the time discretization
by fixing $h=1/128$ and  $T=0.01$. The numerical results
with different time step $\Delta t$ are presented in Table~4.5. We observe that the scheme
(\ref{equ:simplescheme}) has a first-order accuracy in time discretization.

\begin{table}[!htp]
\vskip -0.5cm
\small{\caption{\emph{Spatial errors and convergence rates by scheme (\ref{equ:simplescheme}) for Example 2}}}
\begin{center}
{\small
\begin{tabular}{l|ll|ll|ll|ll}
\hline
$\quad h$&$\quad \quad e_0$   &~$r$  &$\quad\quad e_1$    &~$r$   &$\quad\quad e_{1,r}$ &~$r$ &$\quad\quad e_{2}$ &~$r$ \\
\hline
 $~1/16$ &$1.64\times10^{-2}$ &      &$5.23\times10^{-1}$ &        &$1.31\times10^{-1}$  &    &$4.94\times10^{-0}$  &\\
 $~1/32$ &$3.59\times10^{-3}$ &$2.2$  &$2.18\times10^{-1}$ &$1.3$  &$2.62\times10^{-2}$  &2.3 &$1.08\times10^{-0}$  &2.2\\
 $~1/64$ &$8.71\times10^{-4}$ &$2.0$  &$1.05\times10^{-1}$ &$1.1$  &$6.55\times10^{-3}$  &2.0 &$3.02\times10^{-1}$  &1.9\\
 $1/128$ &$2.16\times10^{-4}$ &$2.0$  &$5.02\times10^{-2}$ &$1.1$  &$1.64\times10^{-3}$  &2.0 &$9.00\times10^{-2}$  &1.8\\
\hline
\end{tabular}}
\end{center}
\label{tab:4}
\end{table}
%

\begin{table}[!htp]
\vskip -0.5cm
\small{\caption{\emph{Temporal errors  and convergence rate by scheme (\ref{equ:simplescheme}) for Example 1}}}
\begin{center}
{\small
\begin{tabular}{l|llll|l}
\hline
$\Delta t$ &$~~~~10^{-4}$      &$~~10^{-4}/2$  &$~~10^{-4}/2^2$    &$~~10^{-4}/2^3$   &~$r$ \\
\hline
 $e_0$ &$ 1.80\times10^{-4}$ &$9.32\times10^{-5}$       &$4.72\times10^{-5}$ &  $2.43\times10^{-5}$       &1.0\\
\hline
\end{tabular}}
\end{center}
\end{table}

\vskip -0.5cm
\subsection{Spinodal decomposition}
In this subsection, we numerically solve the Cahn-Hilliard equation to show the spinodal decomposition:
a process or phenomenon to rapid unmix a mixture of liquids or solids from one thermodynamic phase,
to form two coexisting phases. In the following examples, we apply the scheme \eqref{equ:simplescheme}  on the uniform triangular mesh
with the space stepsize $h=1/128$ and time stepsize $\Delta t=10^{-3}$.
Since the numerical results computed by the scheme \eqref{equ:recnit} on unstructured meshes
are similar to those by \eqref{equ:simplescheme}, they will be not reported here.

{\bf Example 3:} We consider the Cahn-Hilliard equation suggested in \cite{ZW2010} where
the parameter $\epsilon=0.02$ and the initial value is given by
\begin{equation*}
  u_0=10^{-3}\sin^3\frac{\pi x}{4h}\sin^3\frac{\pi y}{4h},~(x,y)\in(0,8h)\times(0,8h).
\end{equation*}

\begin{figure}
   \centering
   \subcaptionbox{\label{fig:tu20}t=0}
  {\includegraphics[width=0.32\textwidth]{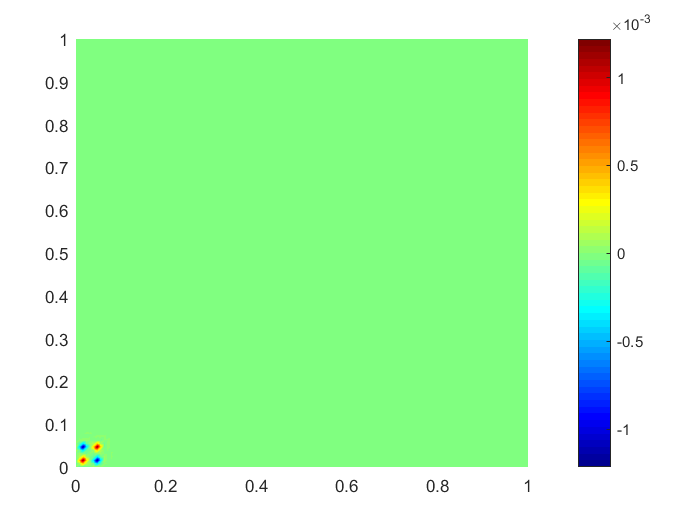}}
  \subcaptionbox{\label{fig:tu21}t=0.01}
   {\includegraphics[width=0.32\textwidth]{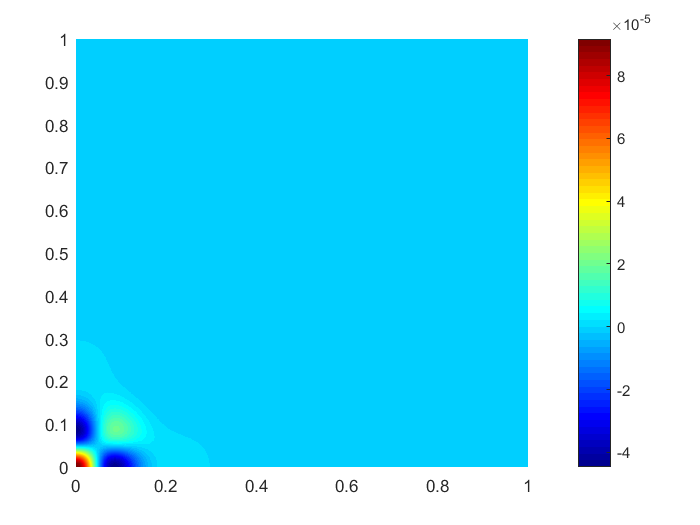}}
  \subcaptionbox{\label{fig:tu22}t=0.1}
  {\includegraphics[width=0.32\textwidth]{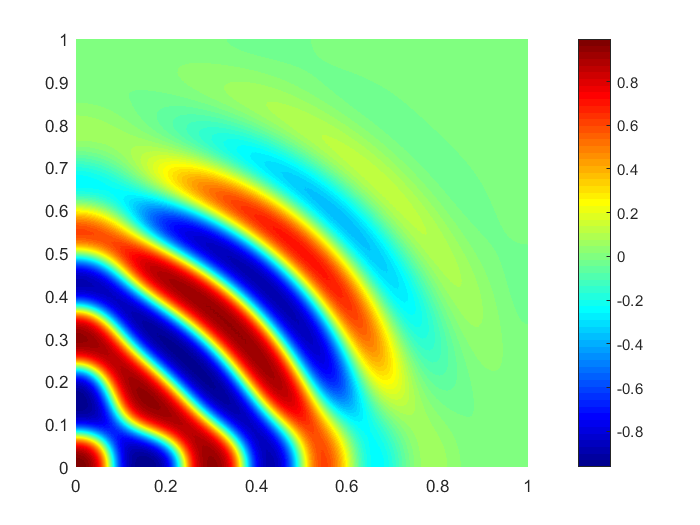}}\\
   \subcaptionbox{\label{fig:tu23}t=0.5}
  {\includegraphics[width=0.32\textwidth]{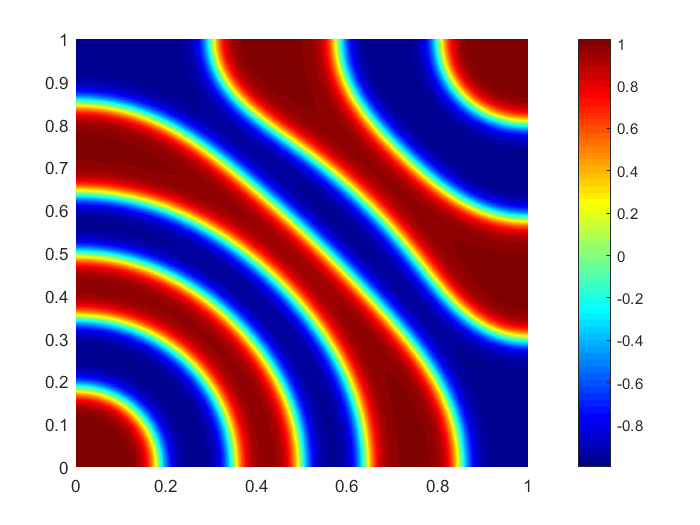}}
  \subcaptionbox{\label{fig:tu24}t=1}
   {\includegraphics[width=0.32\textwidth]{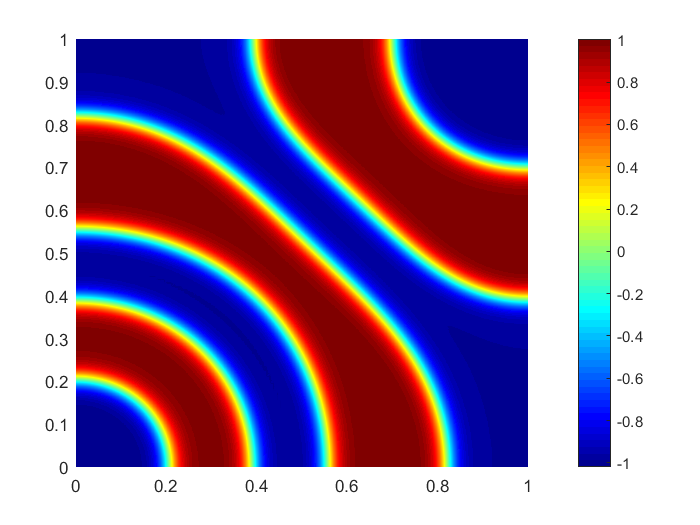}}
  \subcaptionbox{\label{fig:tu25}t=10}
  {\includegraphics[width=0.32\textwidth]{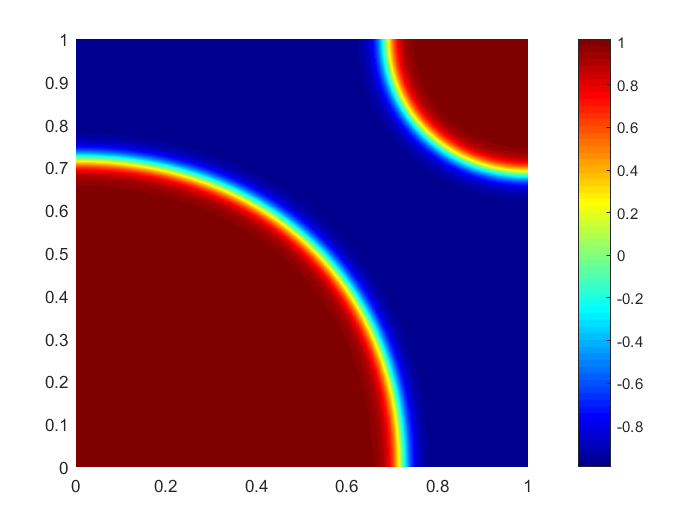}}
  \vskip -0.5cm
   \caption{\small{Example 4, \emph{spinoidal decomposition at six fixed time.}}}
   \label{fig:test3}
  \end{figure}

\begin{figure}[!htp]
   \centering
   \subcaptionbox{\label{fig:ne3}}
  {\includegraphics[width=0.32\textwidth]{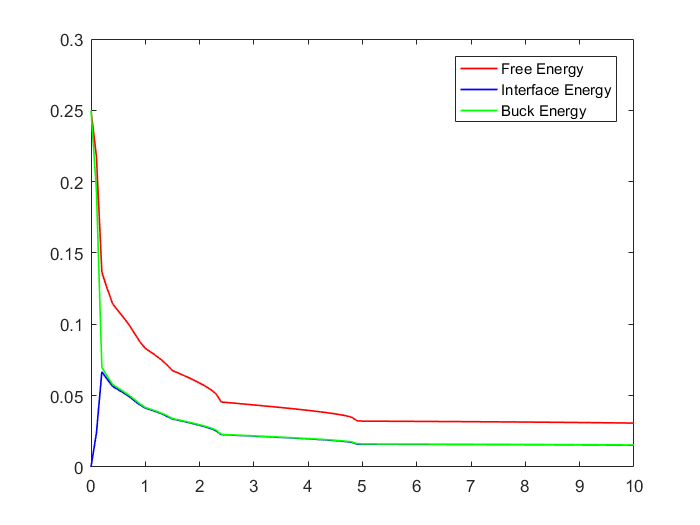}}
  \subcaptionbox{\label{fig:nm3}}
   {\includegraphics[width=0.32\textwidth]{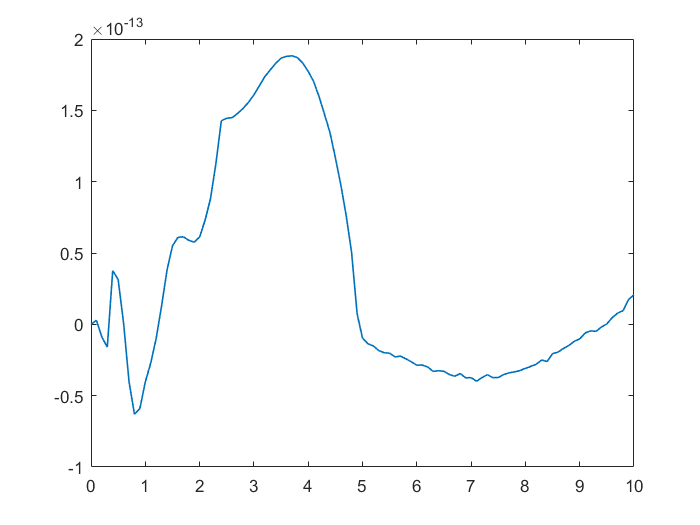}}
   \subcaptionbox{\label{fig:nmax3}}
   {\includegraphics[width=0.32\textwidth]{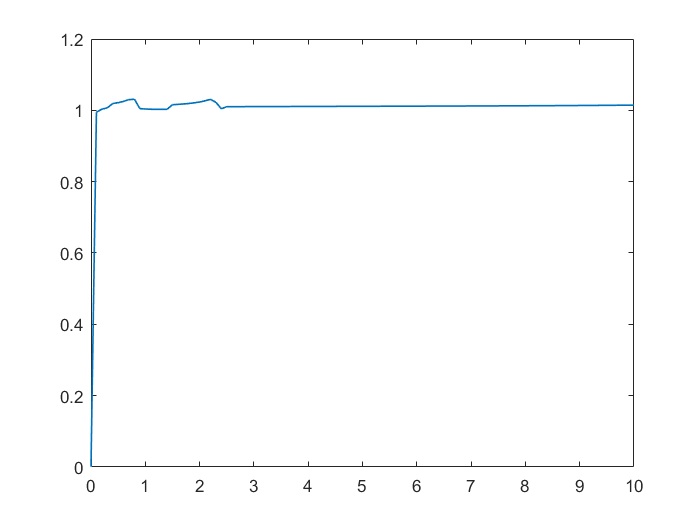}}
   \vskip -0.5cm
   \caption{\small{Example 4 (a): \emph{Energies evolution}, (b): \emph{Mass evolution}, (c): \emph{Evolution of the maximum-norm of the solution.}}}
\end{figure}

We depict the phases at six different times in Fig. \ref{fig:test3}. Note that
typical phase transition phenomena can be clearly observed from these pictures. Moreover, we find that
under a small perturbation($u_0$ is small near the origin), the spinodal decomposition occurs and then coarsens,
and after a period of evolution, the two coexisting phases become stable.  Note that, compared to the
subsequent motion, the initial separation occurs over a very small time scale.
Moreover, the evolution of the energies, including bulk energy and interfacial energy,  is shown in Fig. \ref{fig:ne3},
 the development of the mass is displayed in Fig. \ref{fig:nm3}, the maximum-norm of the numerical solution is illustrated in Fig. \ref{fig:nmax3}.
 Apparently, the presented method almost preserves the properties of energy dissipation and mass conservation, while the  numerical solution itself is
 uniformly bounded.

%

\begin{figure}
   \centering
   \subcaptionbox{\label{fig:tu10}t=0}
  {\includegraphics[width=0.32\textwidth]{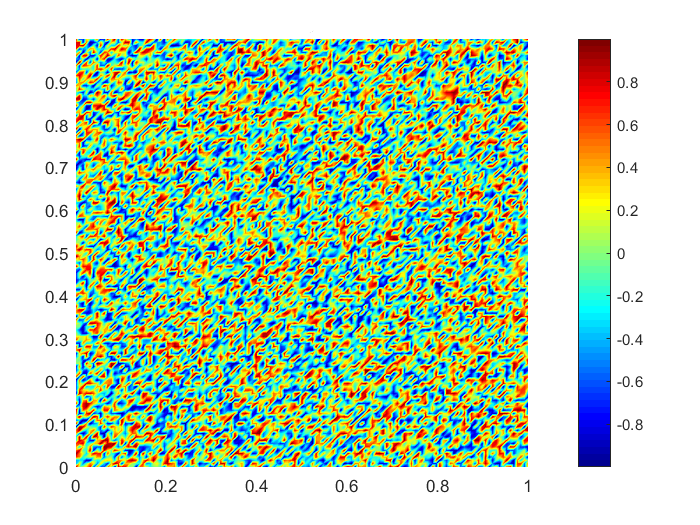}}
  \subcaptionbox{\label{fig:tu11}t=0.01}
   {\includegraphics[width=0.32\textwidth]{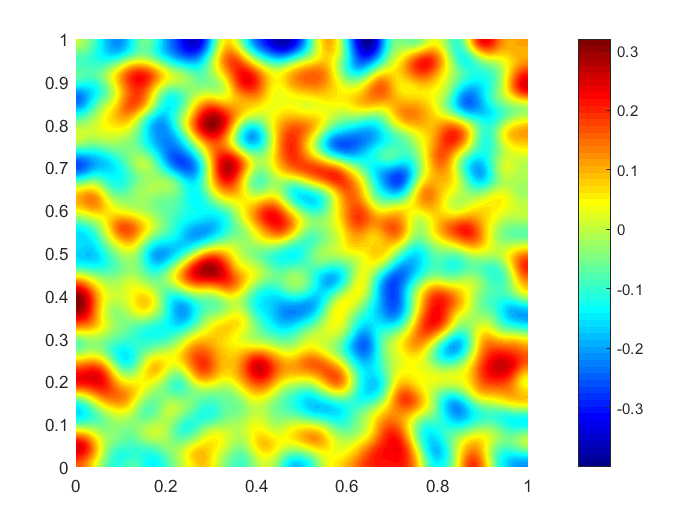}}
  \subcaptionbox{\label{fig:tu12}t=0.1}
  {\includegraphics[width=0.32\textwidth]{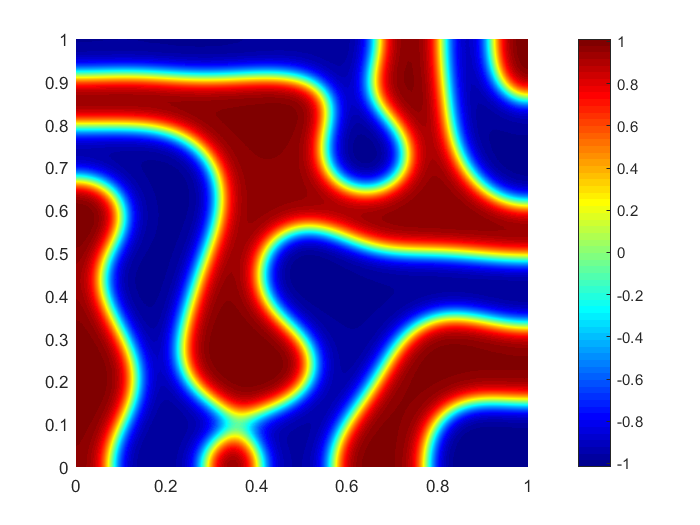}}\\
   \subcaptionbox{\label{fig:tu13}t=0.5}
  {\includegraphics[width=0.32\textwidth]{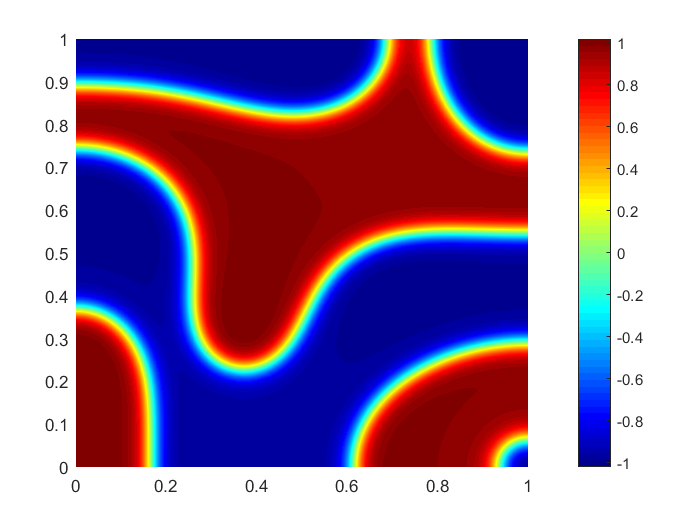}}
  \subcaptionbox{\label{fig:tu14}t=1}
   {\includegraphics[width=0.32\textwidth]{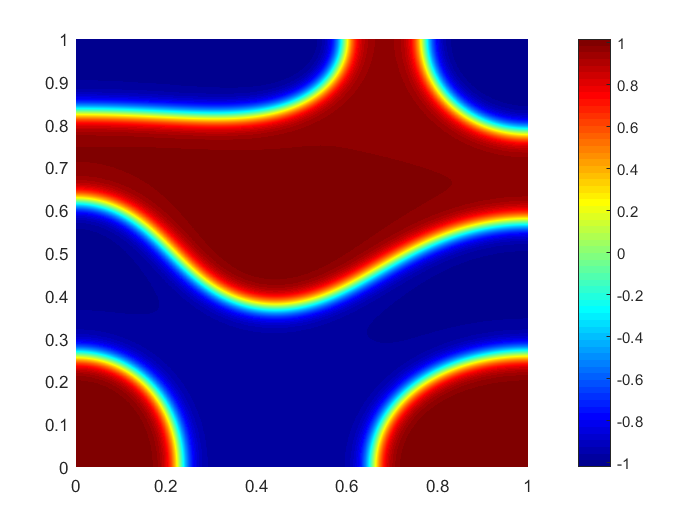}}
  \subcaptionbox{\label{fig:tu15}t=10}
  {\includegraphics[width=0.32\textwidth]{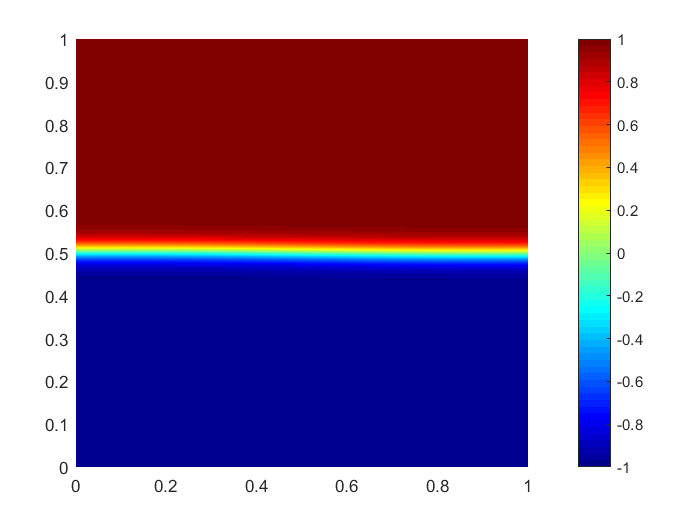}}
  \vskip -0.5cm
   \caption{\small{Example 4, \emph{spinoidal decomposition at six fixed time.}}}
   \label{fig:test4}
  \end{figure}

\begin{figure}[!htp]
   \centering
   \subcaptionbox{\label{fig:ne4}}
  {\includegraphics[width=0.32\textwidth]{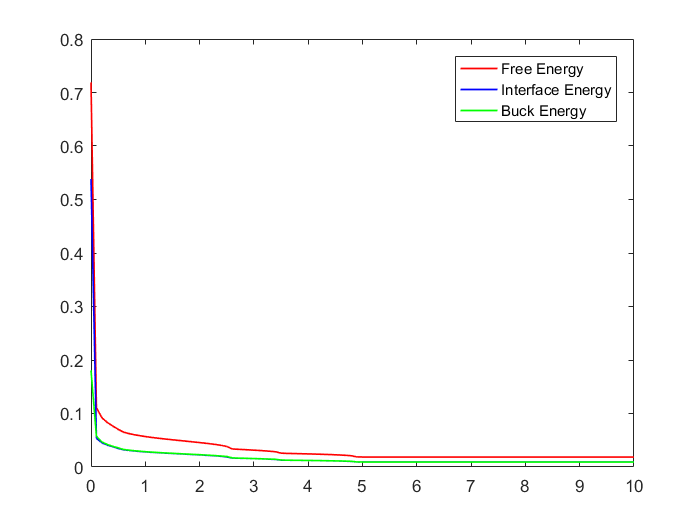}}
  \subcaptionbox{\label{fig:nm4}}
   {\includegraphics[width=0.32\textwidth]{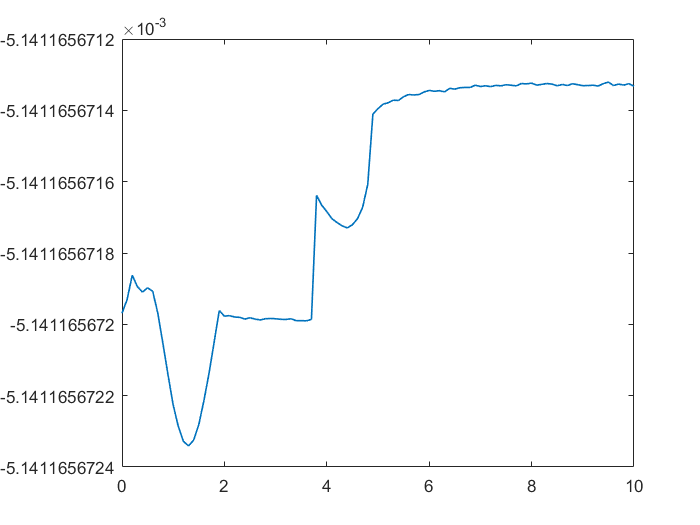}}
   \subcaptionbox{\label{fig:nmax4}}
   {\includegraphics[width=0.32\textwidth]{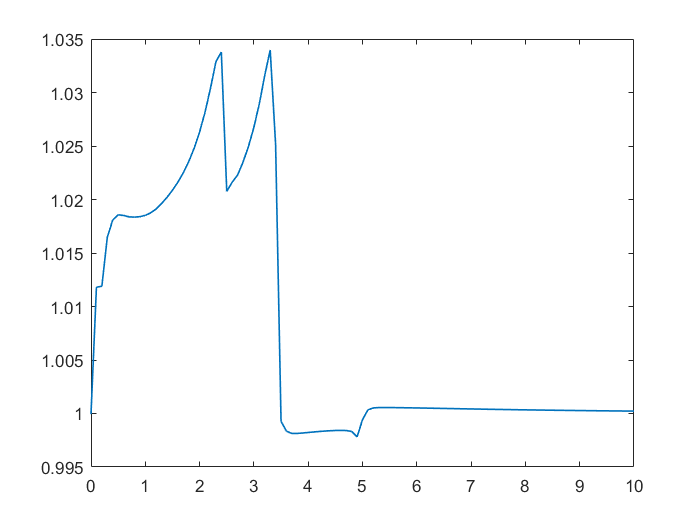}}
   \vskip -0.5cm
   \caption{\small{Example 4 (a): \emph{Energies evolution}, (b): \emph{Mass evolution}, (c): \emph{Evolution of the maximum-norm of the solution.}}}
\end{figure}

{\bf Example 4:} We consider the Cahn-Hilliard equation suggested in \cite{AVSV2016} where
the initial date $u_0$ is a random value field which is uniformly distributed between $-1$ and $1$.
The parameter $\epsilon$ is set to be $0.02$.
We depict the phase evolution of Example 4  in Fig.\ref{fig:test4}.
The process of phase evolution is similar to that in Example 3. That is,  the spinodal decomposition
takes place very early, and after a brief period of evolution, the separation becomes very slow.
Fig. \ref{fig:test4} is also in good agreement with the one presented in \cite{AVSV2016} by using the $C^1$ virtual element method.
The discrete energies and mass are shown in Fig. \ref{fig:ne4} and Fig. \ref{fig:nm4}.
The maximum norm of the approximate solution is displayed in Fig. \ref{fig:nmax4}.
These numerical results reveal that our numerical scheme is energy stable.

%
%

\vskip -0.5cm
\subsection{Evolution of interfaces}
In this subsection, we focus on tracking the evolution of initial data's interfaces, including
a cross-shaped, an elliptic-shaped and two circles-shaped interfaces between phases.  In all examples,
we use the uniform mesh with mesh size $h = 1/128$ and time step size $\Delta t=5\times 10^{-5}$.

{\bf Example 5: } We consider the Cahn-Hilliard equation suggested in \cite{CPMP2016} with $\epsilon=0.01$ and the initial value
\begin{equation*}
u_0(x,y)=
\begin{cases}
~~0.95, & \text{if}~5|(y-0.5)-(x-0.5)|+ |\frac{2}{5}(x-0.5)-(y-0.5)|<1, \\
~~0.95, & \text{if}~5|(x-0.5)-(y-0.5)|+ |\frac{2}{5}(y-0.5)-(x-0.5)|<1, \\
-0.95,  &  \text{otherwise}.
\end{cases}
\end{equation*}

\begin{figure}
   \centering
   \subcaptionbox{\label{fig:tu50}t=0}
  {\includegraphics[width=0.32\textwidth]{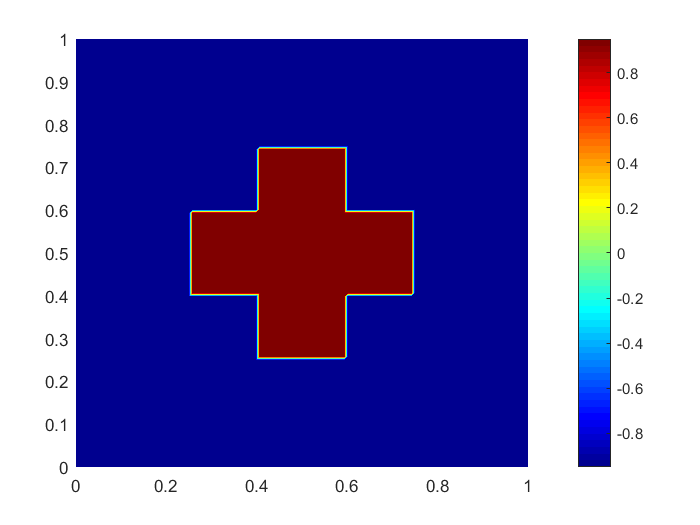}}
  \subcaptionbox{\label{fig:tu51}t=0.005}
   {\includegraphics[width=0.32\textwidth]{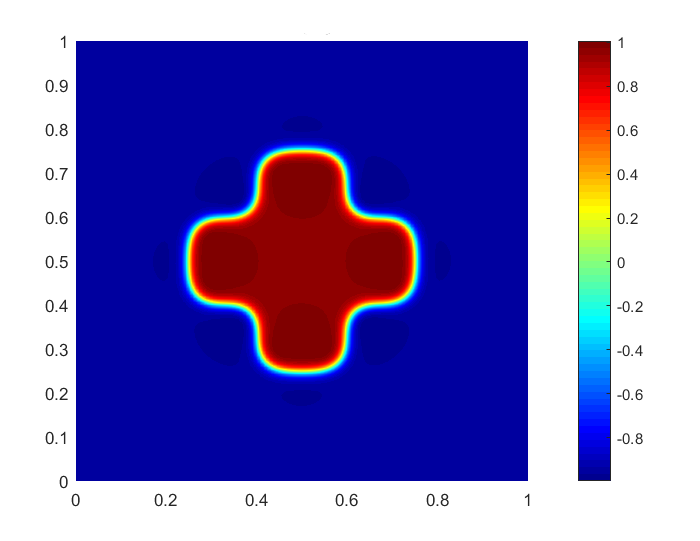}}
  \subcaptionbox{\label{fig:tu52}t=0.01}
  {\includegraphics[width=0.32\textwidth]{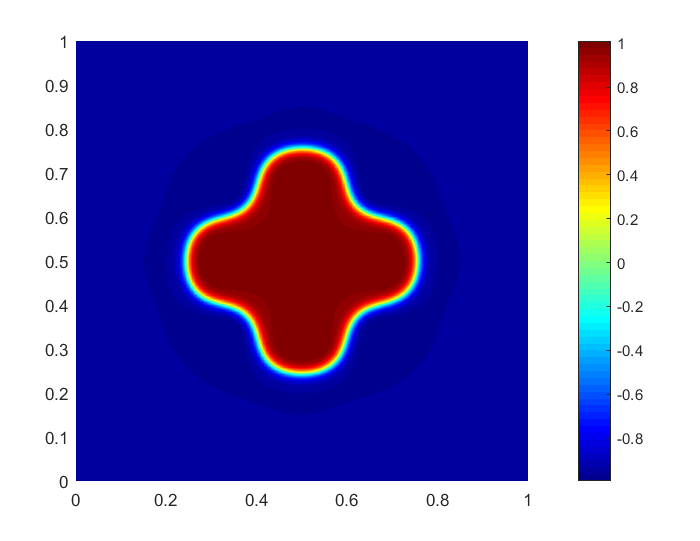}}\\
   \subcaptionbox{\label{fig:tu53}t=0.05}
  {\includegraphics[width=0.32\textwidth]{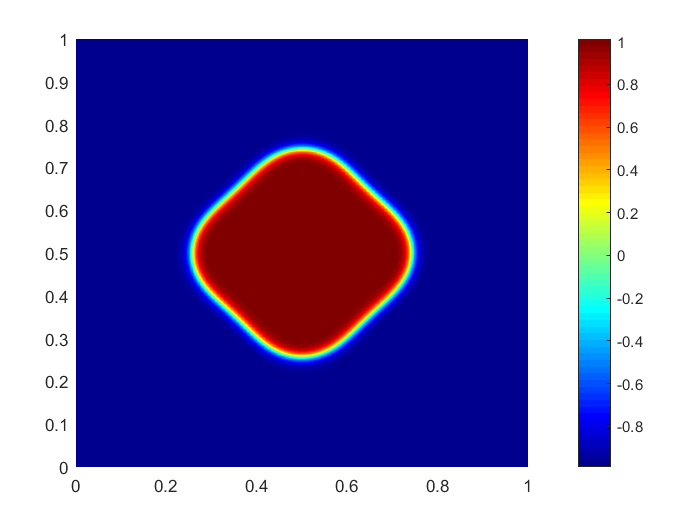}}
  \subcaptionbox{\label{fig:tu54}t=0.1}
   {\includegraphics[width=0.32\textwidth]{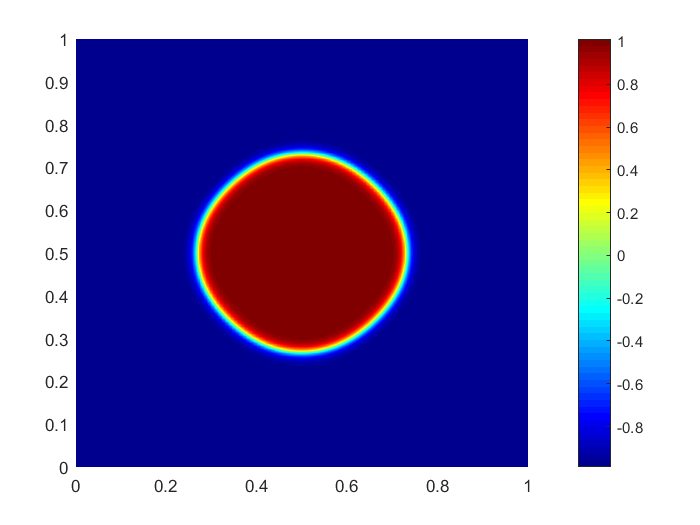}}
  \subcaptionbox{\label{fig:tu55}t=1}
  {\includegraphics[width=0.32\textwidth]{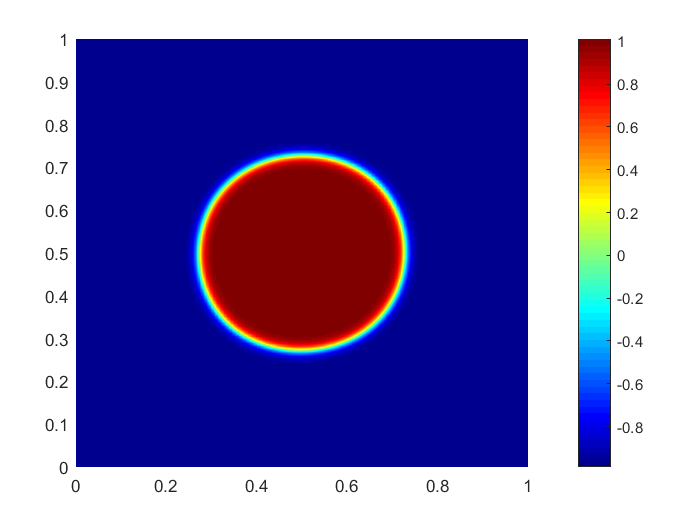}}
  \vskip -0.5cm
   \caption{\small{Example 5 \emph{Evolution of a cross-shaped interface at six temporal frames.}}}
   \label{fig:test5}
\end{figure}

\begin{figure}[!htp]
   \centering
   \subcaptionbox{\label{fig:ne5}}
  {\includegraphics[width=0.32\textwidth]{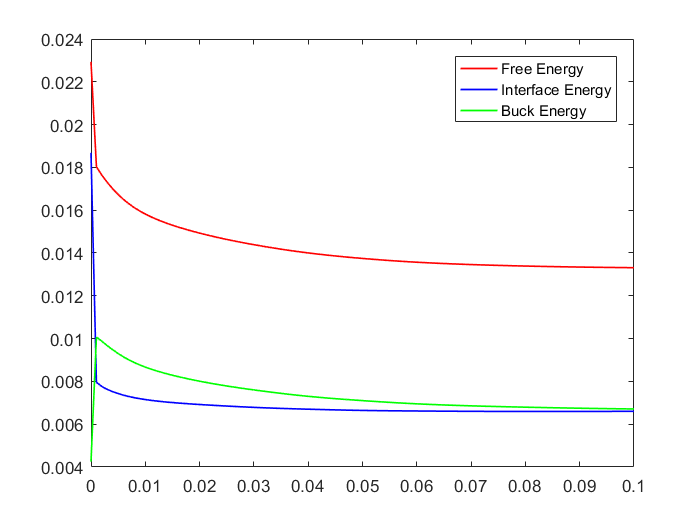}}
  \subcaptionbox{\label{fig:nm5}}
   {\includegraphics[width=0.32\textwidth]{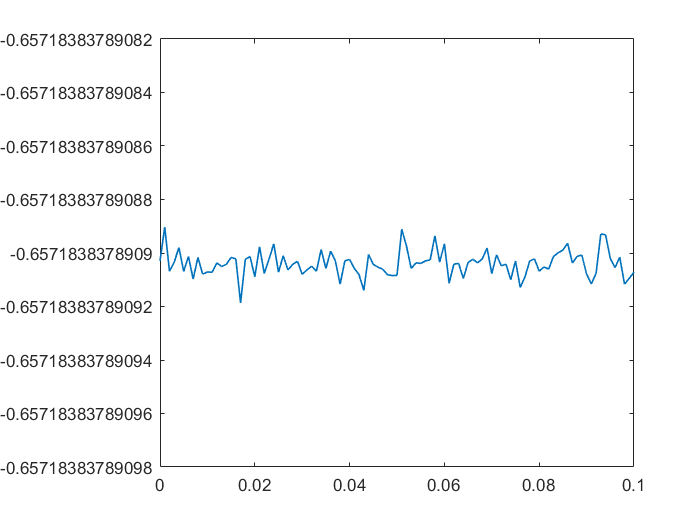}}
   \subcaptionbox{\label{fig:nmax5}}
   {\includegraphics[width=0.32\textwidth]{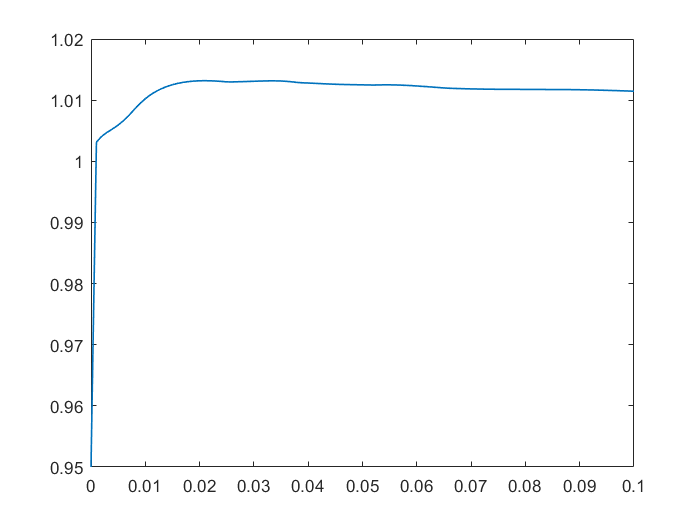}}
   \vskip -0.5cm
    \caption{\small{Example 5 (a) \emph{Energies evolution } (b)\emph{Mass evolution}, (c)\emph{Evolution of the maximum norm of the approximate solution.}}}
   \label{fig:energytest5}
\end{figure}

From Fig. 4.6,  we observe that the cross-shaped interface evolves toward a steady
circular interface. From Fig. \ref{fig:energytest5}, we observe that the
mass is well preserved, the energy is dissipative and the maximum norm of the
approximate solution is controlled. Comparing Fig. \ref{fig:test5} with the
numerical results given in \cite{CPMP2016} computed  by the mixed FEM, they are in good agreement.

{\bf Example 6:} We consider the Cahn-Hilliard equation with a piecewise constant initial data $u_0$ whose
jump set has a shape of an ellipse:
\begin{equation*}
u_0(x,y)=
\begin{cases}
~~0.95, & \text{if}~81(x-0.5)^2+9(y-0.5)^2<1, \\
-0.95,  &  \text{otherwise}.
\end{cases}
\end{equation*}
and the parameter is set to be $\epsilon=0.01$.

\begin{figure}
   \centering
   \subcaptionbox{\label{fig:tu40}t=0}
  {\includegraphics[width=0.32\textwidth]{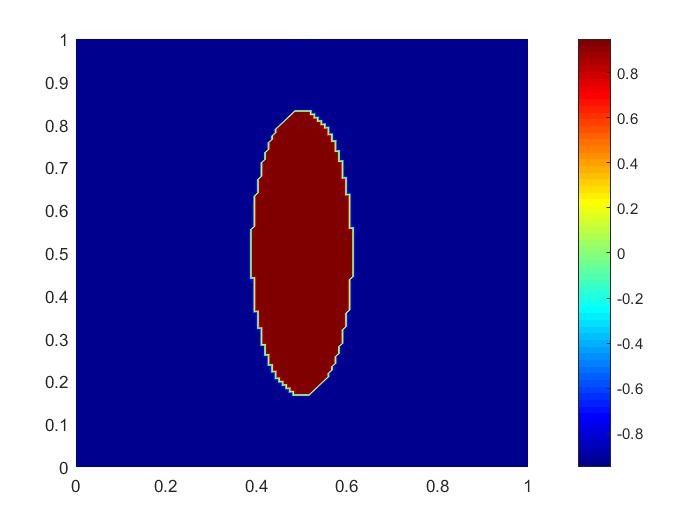}}
  \subcaptionbox{\label{fig:tu41}t=0.003}
   {\includegraphics[width=0.32\textwidth]{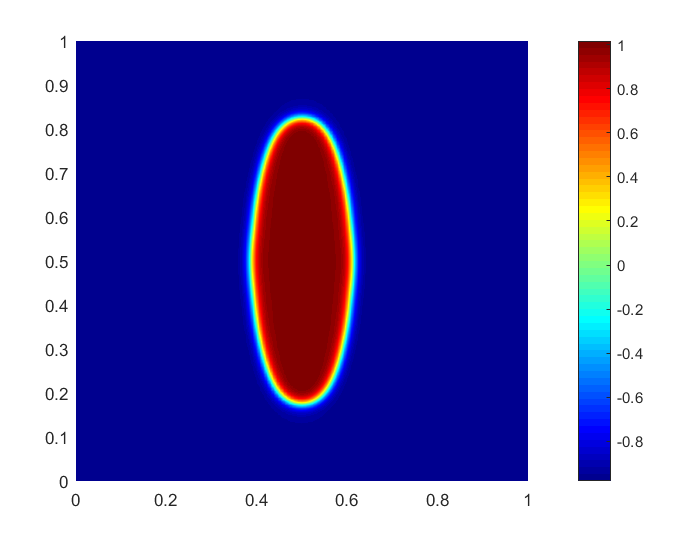}}
  \subcaptionbox{\label{fig:tu42}t=0.05}
  {\includegraphics[width=0.32\textwidth]{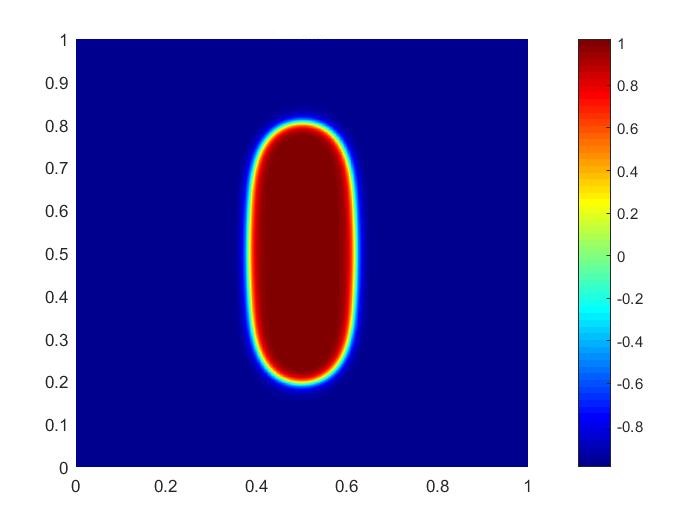}}\\
   \subcaptionbox{\label{fig:tu43}t=0.1}
  {\includegraphics[width=0.32\textwidth]{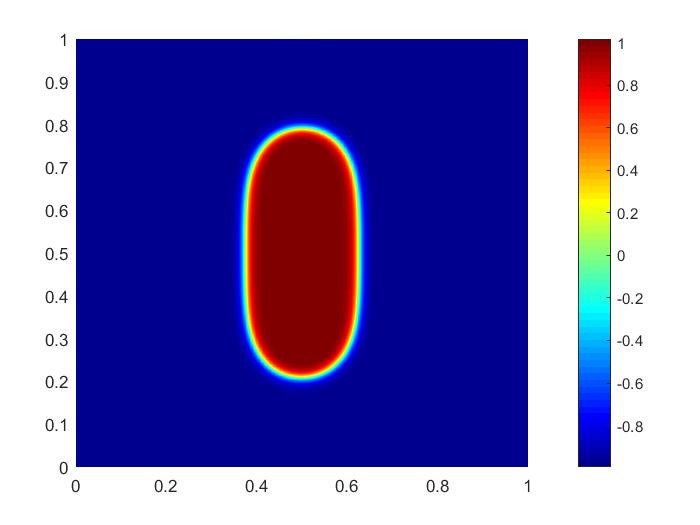}}
  \subcaptionbox{\label{fig:tu44}t=0.3}
   {\includegraphics[width=0.32\textwidth]{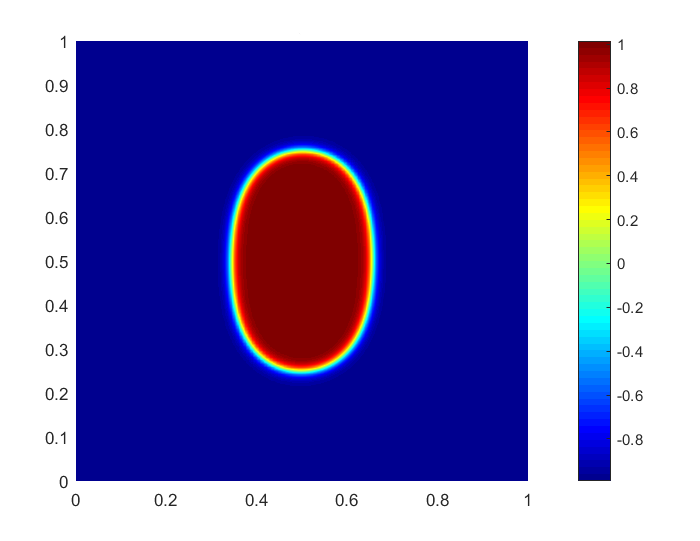}}
  \subcaptionbox{\label{fig:tu45}t=1}
  {\includegraphics[width=0.32\textwidth]{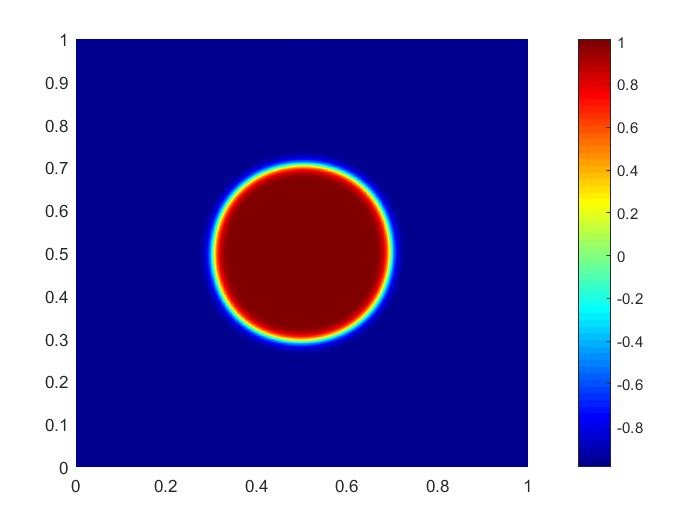}}
   \caption{\small{Example 6 \emph{Evolution of a cross-shaped interface at six temporal frames.}}}
   \label{fig:test6}
\end{figure}

\begin{figure}[!htp]
   \centering
   \subcaptionbox{\label{fig:ne6}}
  {\includegraphics[width=0.32\textwidth]{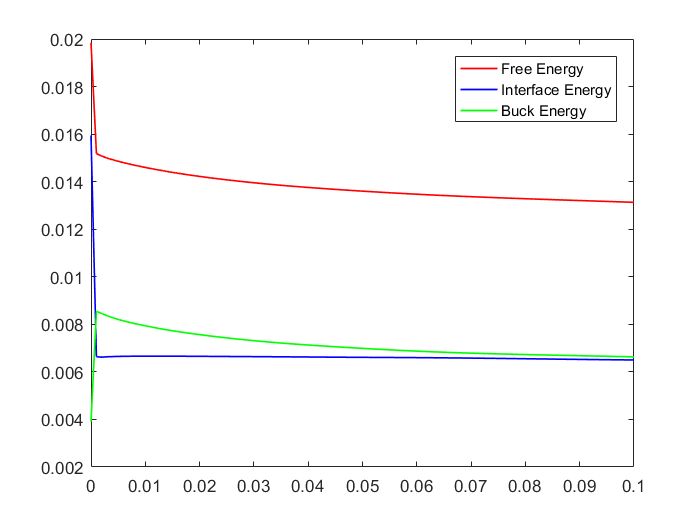}}
  \subcaptionbox{\label{fig:nm6}}
   {\includegraphics[width=0.32\textwidth]{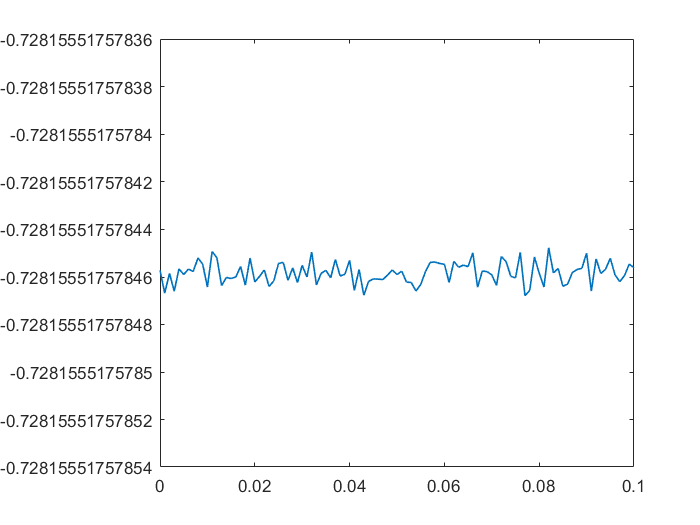}}
   \subcaptionbox{\label{fig:nmax6}}
   {\includegraphics[width=0.32\textwidth]{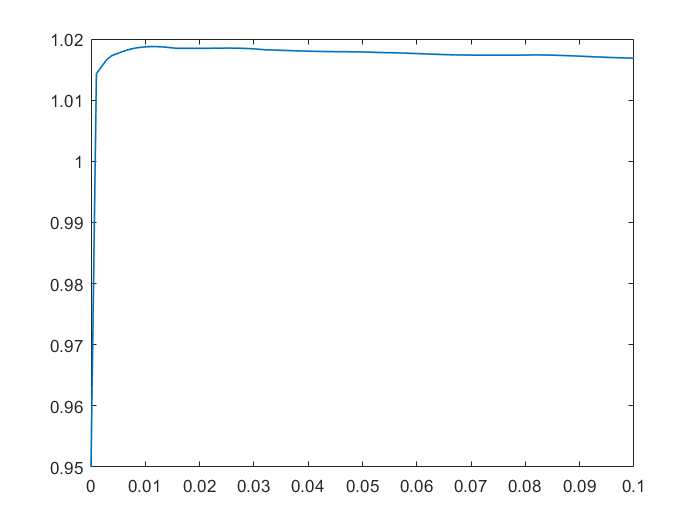}}
   \vskip -0.5cm
    \caption{\small{Example 6 (a) \emph{Energies evolution } (b)\emph{Mass evolution}, (c)\emph{Evolution of the maximum norm of the approximate solution.}}}
      \label{fig:energytest6}
\end{figure}

The numerical results are presented in Fig. \ref{fig:test6}. As in the previous example, we found that
the initial interface
evolves to a steady state exhibiting a circular interface. Moreover, the features of the mass and energies are also captured,
as shown in Fig. \ref{fig:energytest6}.


{\bf Example 7:} We consider the Cahn-Hilliard equation with $\epsilon=0.025$ on the domain $[-1,1]^2$ and  the initial value
\begin{equation*}
  u_0(x)=\tanh (\frac{1}{\sqrt{2}\epsilon} \min\{ \sqrt{(x+0.3)^2+y^2}-0.3,\sqrt{(x-0.3)^2+y^2}-0.25\}).
\end{equation*}

\begin{figure}[!htp]
   \centering
   \subcaptionbox{\label{fig:tu80}t=0}
  {\includegraphics[width=0.32\textwidth]{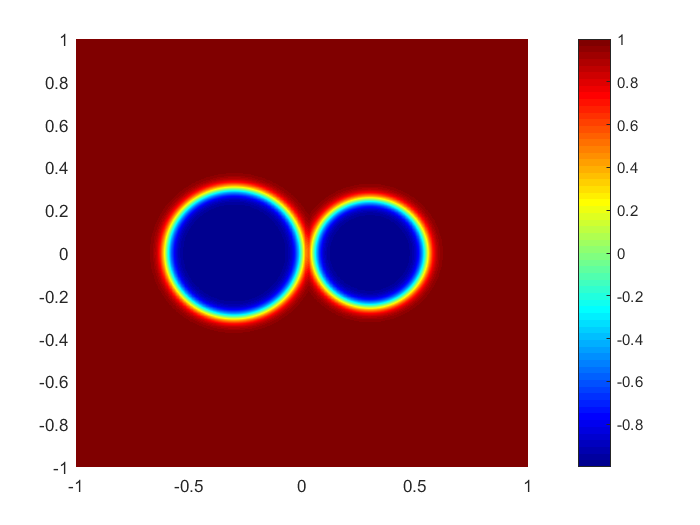}}
  \subcaptionbox{\label{fig:tu81}t=0.001}
   {\includegraphics[width=0.32\textwidth]{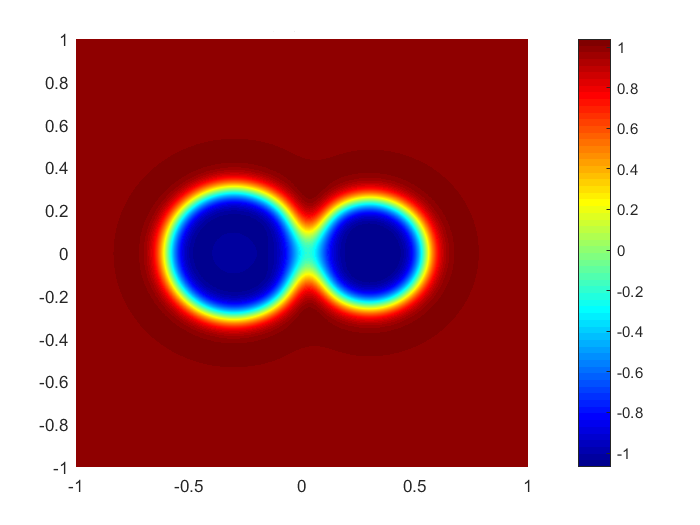}}
  \subcaptionbox{\label{fig:tu82}t=0.005}
  {\includegraphics[width=0.32\textwidth]{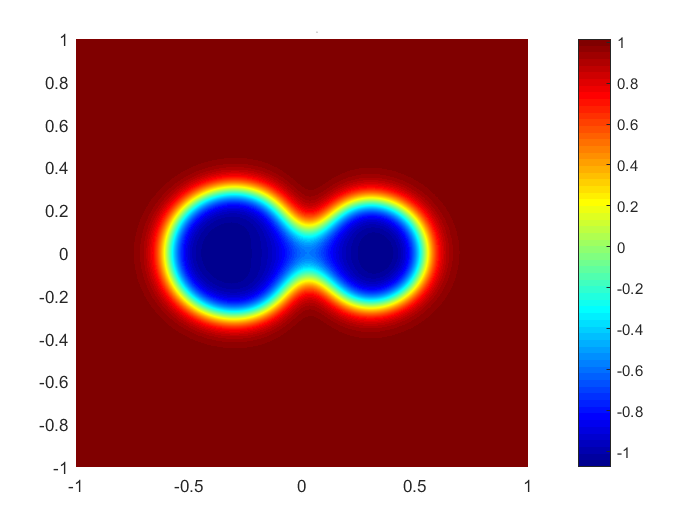}}\\
   \subcaptionbox{\label{fig:tu83}t=0.01}
  {\includegraphics[width=0.32\textwidth]{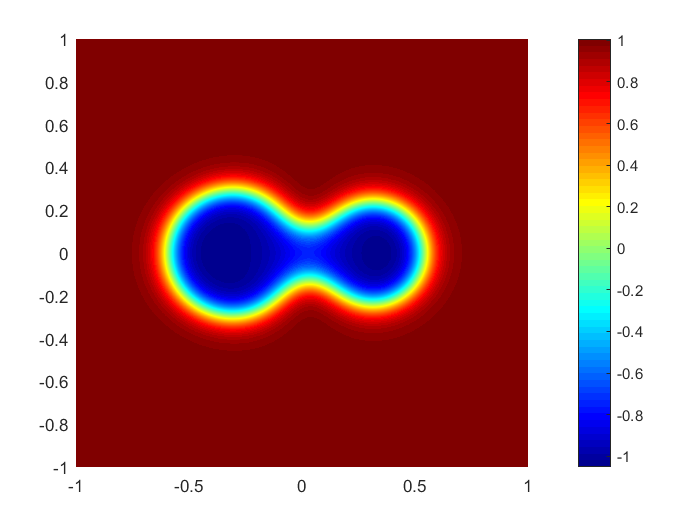}}
  \subcaptionbox{\label{fig:tu84}t=0.05}
   {\includegraphics[width=0.32\textwidth]{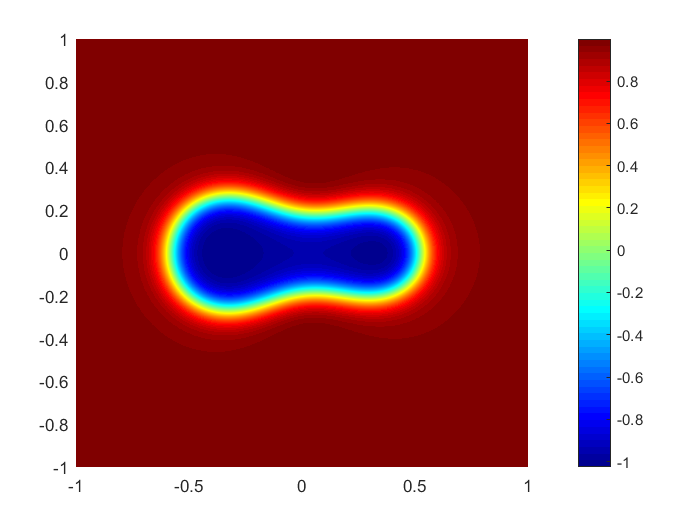}}
  \subcaptionbox{\label{fig:tu85}t=0.1}
  {\includegraphics[width=0.32\textwidth]{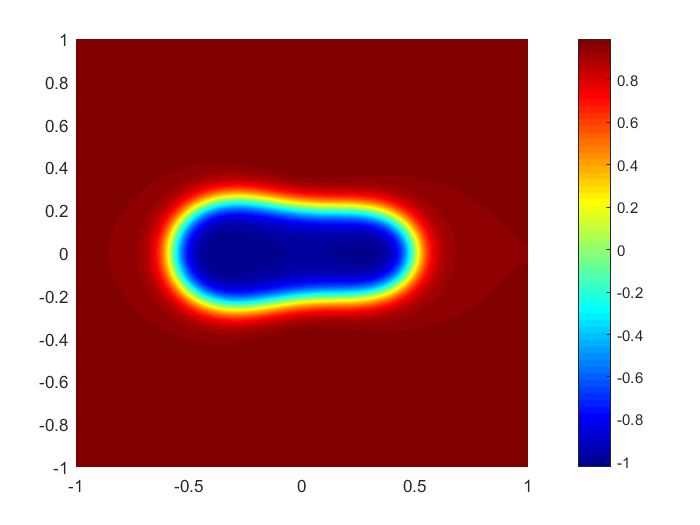}}
   \caption{\small{Example 7 \emph{Evolution of a cross-shaped interface at six temporal frames.}}}
   \label{fig:test7}
\end{figure}

We  generate six snapshots at six fixed time points in Fig. \ref{fig:test7}.
This graph clearly indicates that the two circle interfaces gradually evolve into one circle,
which is consistent with the maximum-norm results obtained in \cite{FLX2016}.  Numerical results depicting  the mass, energies and solution's evolution are
 shown in Fig. 4.11.
\begin{figure}[!htp]
   \centering
   \subcaptionbox{\label{fig:ne7}}
  {\includegraphics[width=0.32\textwidth]{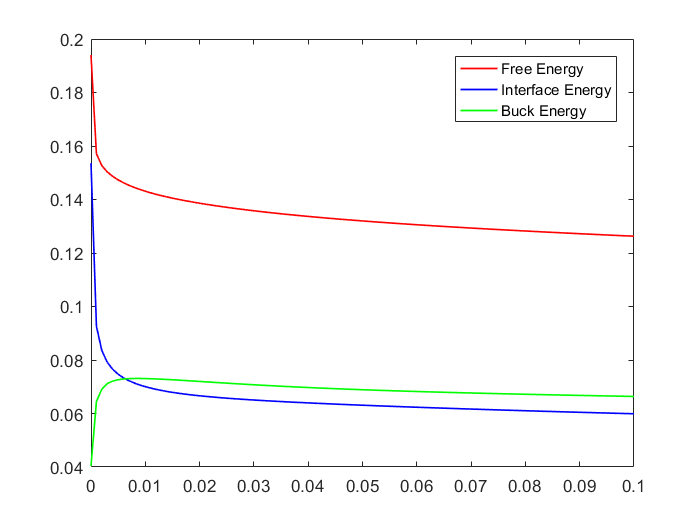}}
  \subcaptionbox{\label{fig:nm7}}
   {\includegraphics[width=0.32\textwidth]{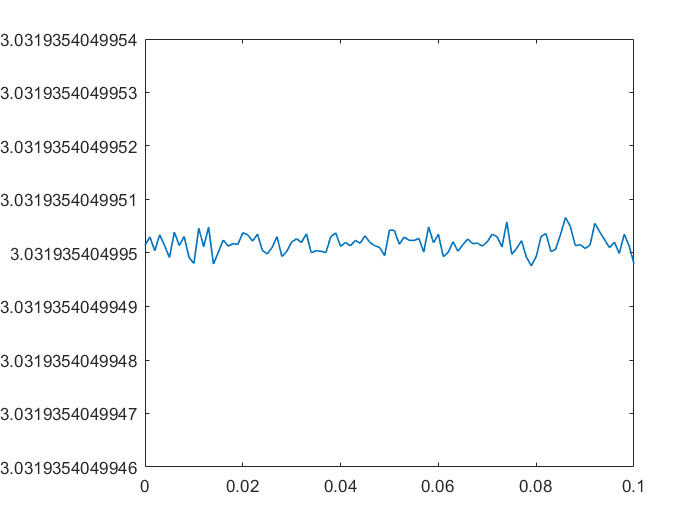}}
   \subcaptionbox{\label{fig:nmax7}}
   {\includegraphics[width=0.32\textwidth]{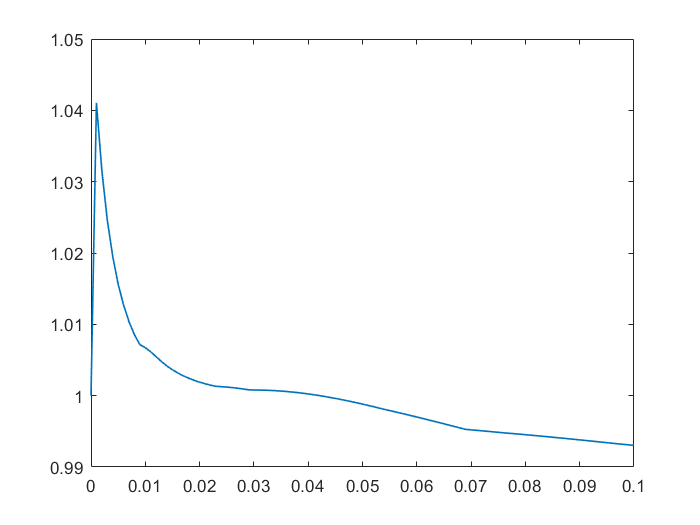}}
    \caption{\small{Example 7 (a) \emph{Energies evolution } (b)\emph{Mass evolution}, (c)\emph{Evolution of the maximum norm of the approximate solution.}}}
\end{figure}

\section{Conclusion}
We designed a $C^0$ linear finite element method
to solve the Cahn-Hilliard equations. This method has a minimum total degree of freedoms and is very simple in implementation.  A series of numerical examples indicate that the new method is stable, efficient and is able to capture some important
physical features such as energy decay and mass conservation during the phase evolution process governed by the Cahn-Hilliard equations.
Meanwhile, the numerical results reveal that  our novel method
has the optimal convergence orders.

Ongoing research topics include a theoretical analysis of the proposed method and
an  extension of the presented method to 3D Cahn-Hilliard
equations and/or other high order differential equations.
For 3D Cahn-Hilliard equation,  we adopt similar  stabilization method \cite{LQ2017b} for time discretization.

\bibliographystyle{siam}
\bibliography{mybibfile}

\begin{thebibliography}{10}

\bibitem{AV2002}
{\sc A.~Agouzal and Y.~Vassilevski}, {\em On a discrete {H}essian recovery for
  {$P_1$} finite elements}, J. Numer. Math., 10 (2002), pp.~1--12.

\bibitem{AVSV2016}
{\sc P.~F. Antonietti, L.~Beir\~ao~da Veiga, S.~Scacchi, and M.~Verani}, {\em A
  {$C^1$} virtual element method for the {C}ahn-{H}illiard equation with
  polygonal meshes}, SIAM J. Numer. Anal., 54 (2016), pp.~34--56.

\bibitem{Ba1972}
{\sc I.~Babu\v{s}ka}, {\em The finite element method with {L}agrangian
  multipliers}, Numer. Math., 20 (1972/73), pp.~179--192.

\bibitem{BCB2003}
{\sc V.~E. Badalassi, H.~D. Ceniceros, and S.~Banerjee}, {\em Computation of
  multiphase systems with phase field models}, J. Comput. Phys., 190 (2003),
  pp.~371--397.

\bibitem{BEG2007}
{\sc Andrea~L. Bertozzi, Selim Esedoḡlu, and Alan Gillette}, {\em Inpainting
  of binary images using the {C}ahn-{H}illiard equation}, IEEE Trans. Image
  Process., 16 (2007), pp.~285--291.

\bibitem{BS2008}
{\sc S.~C. Brenner and L.~R. Scott}, {\em The mathematical theory of finite
  element methods}, vol.~15 of Texts in Applied Mathematics, Springer, New
  York, third~ed., 2008.

\bibitem{CH1958}
{\sc John~W Cahn and John~E Hilliard}, {\em Free energy of a nonuniform system.
  i. interfacial free energy}, The Journal of Chemical Physics, 28 (1958),
  p.~258.

\bibitem{CR2007}
{\sc H.~D. Ceniceros and A.~M. Roma}, {\em A nonstiff, adaptive mesh
  refinement-based method for the {C}ahn-{H}illiard equation}, J. Comput.
  Phys., 225 (2007), pp.~1849--1862.

\bibitem{CPMP2016}
{\sc F.~Chave, D.~A. Di~Pietro, F.~Marche, and F.~Pigeonneau}, {\em A hybrid
  high-order method for the {C}ahn-{H}illiard problem in mixed form}, SIAM J.
  Numer. Anal., 54 (2016), pp.~1873--1898.

\bibitem{CGZZ2017}
{\sc H.~Chen, H.~Guo, Z.~Zhang, and Q.~Zou}, {\em A {$C^0$} linear finite
  element method for two fourth-order eigenvalue problems}, IMA J. Numer.
  Anal., 37 (2017), pp.~2120--2138.

\bibitem{CS2016}
{\sc Y.~Chen and J.~Shen}, {\em Efficient, adaptive energy stable schemes for
  the incompressible {C}ahn-{H}illiard {N}avier-{S}tokes phase-field models},
  J. Comput. Phys., 308 (2016), pp.~40--56.

\bibitem{Ci2002}
{\sc P.~G. Ciarlet}, {\em The finite element method for elliptic problems},
  vol.~40 of Classics in Applied Mathematics, Society for Industrial and
  Applied Mathematics (SIAM), Philadelphia, PA, 2002.
\newblock Reprint of the 1978 original [North-Holland, Amsterdam; MR0520174 (58
  \#25001)].

\bibitem{CR1974}
{\sc P.~G. Ciarlet and P.-A. Raviart}, {\em A mixed finite element method for
  the biharmonic equation},  (1974), pp.~125--145. Publication No. 33.

\bibitem{Co1942}
{\sc R.~Courant}, {\em Variational methods for the solution of problems of
  equilibrium and vibrations}, Bull. Amer. Math. Soc., 49 (1943), pp.~1--23.

\bibitem{DN1991}
{\sc Q.~Du and R.~A. Nicolaides}, {\em Numerical analysis of a continuum model
  of phase transition}, SIAM J. Numer. Anal., 28 (1991), pp.~1310--1322.

\bibitem{EF1989}
{\sc C.~M. Elliott and D.~A. French}, {\em A nonconforming finite-element
  method for the two-dimensional {C}ahn-{H}illiard equation}, SIAM J. Numer.
  Anal., 26 (1989), pp.~884--903.

\bibitem{EFM1989}
{\sc C.~M. Elliott, D.~A. French, and F.~A. Milner}, {\em A second order
  splitting method for the {C}ahn-{H}illiard equation}, Numer. Math., 54
  (1989), pp.~575--590.

\bibitem{EL1992}
{\sc C.~M. Elliott and S.~Larsson}, {\em Error estimates with smooth and
  nonsmooth data for a finite element method for the {C}ahn-{H}illiard
  equation}, Math. Comp., 58 (1992), pp.~603--630, S33--S36.

\bibitem{EZ1986}
{\sc C.~M. Elliott and S.~Zheng}, {\em On the {C}ahn-{H}illiard equation},
  Arch. Rational Mech. Anal., 96 (1986), pp.~339--357.

\bibitem{EGHLMT2002}
{\sc G.~Engel, K.~Garikipati, T.~J.~R. Hughes, M.~G. Larson, L.~Mazzei, and
  R.~L. Taylor}, {\em Continuous/discontinuous finite element approximations of
  fourth-order elliptic problems in structural and continuum mechanics with
  applications to thin beams and plates, and strain gradient elasticity},
  Comput. Methods Appl. Mech. Engrg., 191 (2002), pp.~3669--3750.

\bibitem{FLX2016}
{\sc X.~Feng, Y.~Li, and Y.~Xing}, {\em Analysis of mixed interior penalty
  discontinuous {G}alerkin methods for the {C}ahn-{H}illiard equation and the
  {H}ele-{S}haw flow}, SIAM J. Numer. Anal., 54 (2016), pp.~825--847.

\bibitem{FP2004}
{\sc X.~Feng and A.~Prohl}, {\em Error analysis of a mixed finite element
  method for the {C}ahn-{H}illiard equation}, Numer. Math., 99 (2004),
  pp.~47--84.

\bibitem{GZZ2017b}
{\sc H.~Guo, Z.~Zhang, and R.~Zhao}, {\em Hessian recovery for finite element
  methods}, Math. Comp., 86 (2017), pp.~1671--1692.

\bibitem{GZZZ2016}
{\sc H.~Guo, Z.~Zhang, R.~Zhao, and Q.~Zou}, {\em Polynomial preserving
  recovery on boundary}, J. Comput. Appl. Math., 307 (2016), pp.~119--133.

\bibitem{GZZ2018}
{\sc H.~Guo, Z.~Zhang, and Q.~Zou}, {\em A {$C^0$} linear finite element method
  for biharmonic problems}, J. Sci. Comput., 74 (2018), pp.~1397--1422.

\bibitem{GZZ2018b}
\leavevmode\vrule height 2pt depth -1.6pt width 23pt, {\em A {$C^0$} linear
  finite element method for sixth order elliptic equations}, 2018.
\newblock arXiv:1804.03793v2.

\bibitem{HLT2007}
{\sc Y.~He, Y.~Liu, and T.~Tang}, {\em On large time-stepping methods for the
  {C}ahn-{H}illiard equation}, Appl. Numer. Math., 57 (2007), pp.~616--628.

\bibitem{L2014}
{\sc B.~P. Lamichhane}, {\em A finite element method for a biharmonic equation
  based on gradient recovery operators}, BIT, 54 (2014), pp.~469--484.

\bibitem{Le2007}
{\sc R.~J. LeVeque}, {\em Finite difference methods for ordinary and partial
  differential equations}, Society for Industrial and Applied Mathematics
  (SIAM), Philadelphia, PA, 2007.
\newblock Steady-state and time-dependent problems.

\bibitem{LQ2017}
{\sc D.~Li and Z.~Qiao}, {\em On second order semi-implicit {F}ourier spectral
  methods for 2{D} {C}ahn-{H}illiard equations}, J. Sci. Comput., 70 (2017),
  pp.~301--341.

\bibitem{LQ2017b}
\leavevmode\vrule height 2pt depth -1.6pt width 23pt, {\em On the stabilization
  size of semi-implicit {F}ourier-spectral methods for 3{D} {C}ahn-{H}illiard
  equations}, Commun. Math. Sci., 15 (2017), pp.~1489--1506.

\bibitem{LQT2016}
{\sc D.~Li, Z.~Qiao, and T.~Tang}, {\em Characterizing the stabilization size
  for semi-implicit {F}ourier-spectral method to phase field equations}, SIAM
  J. Numer. Anal., 54 (2016), pp.~1653--1681.

\bibitem{Mo1968}
{\sc L.~S.~D MORLEY}, {\em The triangular equilibrium element in the solution
  of plate bending problems}, Aerosp. Quart., {19} ({1968}), pp.~{149--169}.

\bibitem{easymesh}
{\sc B.~Niceno}, {\em {EasyMesh}: A two-dimensional quality mesh generator}.
\newblock \url{http://web.mit.edu/easymesh_v1.4/www/easymesh.html}, 2001.

\bibitem{Ni1971}
{\sc J.~Nitsche}, {\em \"uber ein {V}ariationsprinzip zur {L}\"osung von
  {D}irichlet-{P}roblemen bei {V}erwendung von {T}eilr\"aumen, die keinen
  {R}andbedingungen unterworfen sind}, Abh. Math. Sem. Univ. Hamburg, 36
  (1971), pp.~9--15.
\newblock Collection of articles dedicated to Lothar Collatz on his sixtieth
  birthday.

\bibitem{OHP2010}
{\sc J.~T. Oden, A.~Hawkins, and S.~Prudhomme}, {\em General diffuse-interface
  theories and an approach to predictive tumor growth modeling}, Math. Models
  Methods Appl. Sci., 20 (2010), pp.~477--517.

\bibitem{PABG2011}
{\sc M.~Picasso, F.~Alauzet, H.~Borouchaki, and P.-L. George}, {\em A numerical
  study of some {H}essian recovery techniques on isotropic and anisotropic
  meshes}, SIAM J. Sci. Comput., 33 (2011), pp.~1058--1076.

\bibitem{SY2010}
{\sc J.~Shen and X.~Yang}, {\em Numerical approximations of {A}llen-{C}ahn and
  {C}ahn-{H}illiard equations}, Discrete Contin. Dyn. Syst., 28 (2010),
  pp.~1669--1691.

\bibitem{VMDDG2007}
{\sc M.-G. Vallet, C.-M. Manole, J.~Dompierre, S.~Dufour, and F.~Guibault},
  {\em Numerical comparison of some {H}essian recovery techniques}, Internat.
  J. Numer. Methods Engrg., 72 (2007), pp.~987--1007.

\bibitem{WKG2006}
{\sc G.~N. Wells, E.~Kuhl, and K.~Garikipati}, {\em A discontinuous {G}alerkin
  method for the {C}ahn-{H}illiard equation}, J. Comput. Phys., 218 (2006),
  pp.~860--877.

\bibitem{WLFC2008}
{\sc S.~M. Wise, J.~S. Lowengrub, H.~B. Frieboes, and V.~Cristini}, {\em
  Three-dimensional multispecies nonlinear tumor growth---{I}: {M}odel and
  numerical method}, J. Theoret. Biol., 253 (2008), pp.~524--543.

\bibitem{XXS2007}
{\sc Y.~Xia, Y.~Xu, and C.-W. Shu}, {\em Local discontinuous {G}alerkin methods
  for the {C}ahn-{H}illiard type equations}, J. Comput. Phys., 227 (2007),
  pp.~472--491.

\bibitem{ZW2010}
{\sc S.~Zhang and M.~Wang}, {\em A nonconforming finite element method for the
  {C}ahn-{H}illiard equation}, J. Comput. Phys., 229 (2010), pp.~7361--7372.

\bibitem{ZN2005}
{\sc Z.~Zhang and A.~Naga}, {\em A new finite element gradient recovery method:
  superconvergence property}, SIAM J. Sci. Comput., 26 (2005), pp.~1192--1213
  (electronic).

\bibitem{ZTZ2013}
{\sc O.~C. Zienkiewicz, R.~L. Taylor, and J.~Z. Zhu}, {\em The finite element
  method: its basis and fundamentals}, Elsevier/Butterworth Heinemann,
  Amsterdam, seventh~ed., 2013.

\bibitem{ZZ1992a}
{\sc O.~C. Zienkiewicz and J.~Z. Zhu}, {\em The superconvergent patch recovery
  and a posteriori error estimates. {I}. {T}he recovery technique}, Internat.
  J. Numer. Methods Engrg., 33 (1992), pp.~1331--1364.

\end{thebibliography}
\end{document}